\documentclass[12pt,a4paper]{article}
\usepackage{
   amsmath,
   amsfonts,
   amssymb,
   graphicx,
   epsfig,
   stmaryrd,
   epstopdf,
   booktabs,
   algorithm,
   algorithmic,
   rotating,
   caption,
   subcaption,
   pbox,
   multirow,
   verbatim,
   pgfplots,
   url,
   setspace,
   tikz
   }
\usepackage{rotating}
\usepackage{hvfloat}
\usepackage[width=1.0\textwidth]{caption}

\pgfplotsset{every axis/.append style={
                    label style={font=\footnotesize},
                    tick label style={font=\tiny}  
                    }}

\def\p{{FIXB}}
\def\pp{\bar{p}}

\newcommand{\ncoppa}[2]{\genfrac{}{}{0pt}{1}{#1}{#2}}

\newtheorem{theorem}{Theorem}

\newtheorem{lemma}[theorem]{Lemma}

        {\hspace*{\fill}$\Box$\par\vspace{4mm}}

\newenvironment{proof}{\noindent{\em Proof.}}%
        {\hspace*{\fill}$\Box$\par\vspace{4mm}}

\begin{document}

\title{Flow Shop Scheduling with Inter-Stage Flexibility \\and Blocking Constraints}
\author{Gaia Nicosia\thanks{Dipartimento di Ingegneria   Civile, Informatica e delle Tecnologie Aeronautiche,
Universit\`a degli Studi ``Roma Tre'', via della Vasca Navale 79,
00146 Rome, Italy, e-mail: gaia.nicosia@uniroma3.it}
\and  
Andrea Pacifici\thanks{Dipartimento di Ingegneria Civile e Ingegneria
Informatica, Universit\`a degli Studi di Roma ``Tor Vergata'', Via del Politecnico 1, 00133 Rome, Italy, e-mail:
andrea.pacifici@uniroma2.it, anna.russo.russo@uniroma2.it, cecilia.salvatore@uniroma2.it}
\and
Ulrich Pferschy\thanks{Department of Operations and Information Systems, University of Graz, Universitaetsstrasse
15, 8010 Graz, Austria. e-mail: ulrich.pferschy@uni-graz.at}
\and
Anna Russo Russo\footnotemark[2] 
\and
Cecilia Salvatore\footnotemark[2]
}

\date{}
\maketitle

\abstract{We investigate a scheduling problem arising from a material handling and processing problem in a production line of an Austrian company building prefabricated house walls.
The addressed problem is a permutation flow shop with blocking constraints in which  the
machine of at least one stage can process a number operations of two other stages in the system. 
This situation is usually referred to as multi-task or inter-stage flexibility.
The problem is in general NP-hard, but we derive a number of special cases that can be solved in polynomial time. 
For the general case, we present a variety of heuristic algorithms, with a focus on matheuristics that are grounded in two different mixed-integer linear programming (MIP) formulations of the problem. 
These matheuristics leverage the strengths of exact optimization techniques while introducing flexibility to address limits on computation time.
To assess the performance of the proposed approaches, we conduct an extensive computational study on randomly generated test cases based on real-world instances.
}


\section{Introduction}

The flow-shop scheduling problem, in its basic formulation, consists of determining how to process $n$ jobs on $m$ machines placed along a line in order to minimize the makespan. 
Each job consists of exactly $m$ operations, which must be executed on their respective machines, in the strict order imposed by the line. Once started, the processing of an operation cannot be interrupted and no machine can perform more than one operation at a time.  
Flow shop can also be seen as a special case of job-shop with a strict order for all operations.
These problems have been studied since the 
Fifties (see for instance~\cite{Johnson54}) due to their importance in various 
industries. In fact, advancements and innovations in flow shop scheduling 
methodologies  significantly influence the performance and competitiveness 
of industrial processes due to their direct impact on operational efficiency, 
production costs, and overall resource utilization. 

A special situation, denoted as  \emph{blocking} constraint, arises in flow shops when a job, that has completed on a machine, blocks it (i.e., makes it not available) until the next downstream machine becomes free.
This typically occurs when there is no buffer between two consecutive machines $M_k$ and $M_{k+1}$ and therefore if $M_{k+1}$ is occupied by a job then the subsequent job $j$ in the sequence blocks $M_k$ even if $j$ already completed its processing on $M_k$.
Since a job cannot overtake another one, all machines process the jobs in the same order. Hence, blocking also implies that the problem becomes a so-called \emph{permutation} flow shop~\cite{blazewicz2007handbook} having the property that the job sequence on the first machine remains the same on all the other machines.

In this paper we  address a problem introduced in \cite{us_ods2023} and inspired by a real-world scenario arising in a production line at a company specializing in constructing prefabricated house walls.
It is a generalization of a flow shop with blocking constraints.
In this problem, a number of operations can be allocated to either 
of two consecutive machines in the flow line, and the operations' processing times vary depending on the designated machine. 
In the literature, this situation is referred to as \emph{multi-task flexibility} or \emph{inter-stage flexibility}.  
In the following, we also describe this situation as \emph{flexible operation-to-machine assignment}. 
Machines able to perform multiple types of operations are usually called \emph{flexible}, while we refer to operations that can be processed by more than one machine as \emph{shiftable}.
As a consequence, in addition to the decision concerning the job sequence, as in the classic permutation flow shop, also the decision on the assignment has to be taken. 

Hereafter, after describing the application scenario that motivated this study (Section~\ref{sec:application}), we outline the contributions of this paper and the structure of the remaining sections.
\begin{itemize}
    \item We provide an overview of the current state-of-the-art literature concerning both the permutation flow shop and the flow shop with multi-task flexibility (Section~\ref{sec:literature}).
    \item We give a formal statement of the addressed problem together with a complexity characterization (Section~\ref{sect:statement}).
    \item We study two special  cases of the problem, namely the 2-jobs and the 2-machines with fixed sequence cases, and show that they can be efficiently solved  (Section \ref{sec:special_cases}). 
    \item We propose a number of \emph{matheuristics} exploiting two different  MIP models (Section~\ref{sec:heur_ilp}) for the general problem. 
    We also present a constructive heuristic, based on the iterative insertion of jobs into a partial schedule until a complete solution is built (Section~\ref{sec:constrheu}). 
    \item We perform an extensive computational campaign based on randomly generated data following from the real-world  application and an extension thereof.
    We present the results of this experimental analysis and  provide a thorough comparison of the proposed algorithms efficiency, solution quality, and scalability across diverse problem sizes and conditions, offering valuable insights into their applicability for real-world decision-making scenarios
    (Section~\ref{sec:results}).
\end{itemize}

\subsection{Motivating real-world case}\label{sec:application}

The optimization model addressed in this work is inspired by a material handling problem arising at a production line  of an Austrian company building prefabricated house walls. 

Currently, the production line is organized in such a way that the house walls (or jobs) pass through five different workstations, namely $M_1, M_2,\ldots, M_5$, in order. 
On each of these workstations, several processing steps or tasks (i.e., operations, such as frame building, socket-holes cutting, wall insulation, etc.) are carried out by the company's workers and, in one \emph{robot-station} by a flexible automated machine able to execute certain operations without the assistance of human operators. 
 Due to the general construction principle of house walls the main steps of production follow the same sequence through the individual stations, although operations for each wall on every machine may differ considerably.
These depend on the individual requirements placed on each wall, e.g., kitchen-bathroom connections for the separating walls, provisions for plumbing and wiring, positions of doors and windows.
Accordingly, ($i$) each wall visits the stations in the same sequence, from $M_1$ to $M_5$ (i.e., machine $M_i$, $i=1,2,\dots,5$, is the $i$-th visited station). 
In particular, the robot-station is the third one in the sequence ($M_3$). Moreover, due to the large size and weight of the bulky walls, no intermediate storage is available between the individual stations. 
As a consequence, ($ii$) a wall cannot move to the next work station if this is busy, thus blocking the upstream machines and possibly causing delay in the overall process. 

Flexibility of the automated machines  
allows a more efficient handling of the operations
since some tasks can be transferred from/to the robot-station $M_3$ to/from the upstream or downstream machines $M_2$ and $M_4$. This way, if for instance machine $M_3$ is still busy with the preceding job while $M_2$ has already completed its tasks on the current job, it is possible for $M_2$ to start performing some additional tasks previously assigned to $M_3$ and save time. 

The company has to decide how to schedule wall elements and manage the execution of tasks on the machines based on a given set of orders, i.e., a given set of walls each with certain specific requirements. 
These specifications imply a given set of tasks that must be performed, in a strict order, on the machines. 
Some of these tasks must be processed by a specific machine while some others can be executed either by machines $M_2$/$M_4$ or the flexible automated machine $M_3$.
Obviously the processing time of the latter tasks may vary depending on whether it is performed by $M_3$ or by a human operator on $M_2$ and $M_4$.

In conclusion, the company has to decide $(i)$ the order in which the walls will be processed and $(ii)$ for each wall, the assignment of its shiftable tasks to a suitable machine on the line such that the makespan, i.e., the completion time of the last job on the last machine, is minimized.

\subsection{Related literature}\label{sec:literature}

When considering a flow line with no buffers, blocking flow shop scheduling problems arise. 
Our problem falls in this area and, due to the above illustrated issues, it may be modeled as a permutation flow shop with blocking constraints and some additional flexibility characteristics. 
Hereafter, we  provide a brief overview of the scientific literature concerning makespan minimization flow shop problems, starting with versions that include blocking constraints and then addressing flexible operation-to-machine assignment.

Flow shop problems with blocking constraints and with the objective of minimizing the makespan are NP-hard in the strong sense as soon as the shop has three machines \cite{Hall1996,papa80}, 
while in the case of two machines, the problem is reduced to a special traveling salesman problem \cite{Reddi1972} that can be efficiently solved using the Gilmore and Gomory algorithm \cite{gg1964}. 
Due to the inherent complexity of the general problem, exact methods are typically employed to solve small instances, whereas heuristics and metaheuristics methods are more commonly utilized for larger instances (see, e.g. \cite{Ozolins2019,Suhami1981,Tasgetiren2017,Wang2011}).
A survey on blocking flow shop scheduling problems can be found in \cite{survey_blocking_2019}.

Closely related to our specific flow shop scheduling problem, is the multi-task flexibility or inter-stage flexibility characteristic, which involves the flexible assignment of operations to machines within a flow line~\cite{azaiez,khorasanian}. 
This flexibility allows certain operations to be assigned to more than one machine in the line. These types of problems have been studied since the nineties \cite{liao,pan} and, in particular,  Pan and Chen~\cite{pan} prove that the makespan minimization problem is already NP-hard  when there are two machines and each job consists of exactly two operations.

Notably, given the complexity inherent in multi-task flexibility flow shop problems, numerous researchers focus on specific problem instances, such as those involving only 2 machines and featuring 2 to 3 operations per job (see, for instance, \cite{azaiez,dellacroce2022,khorasanian,Lin2017}).
The case in which each job consists of three operations and
the first and the last operation must be performed on the first and the second machine, while the second operation  can be performed by either one of the machines has been considered in~\cite{dellacroce2022,Gupta2004,Lin2017}.
In~\cite{Gupta2004} a number of approximations results are presented, while in~\cite{dellacroce2022}, the authors propose some MIP models and exact solution methods based on  a constraint generation approach.
In \cite{Lin2017} the special case in which the sequence of jobs is fixed has been proved to be NP-hard and pseudopolynomial algorithms are proposed. 
In \cite{khorasanian} a two-machine preemptive flow shop scheduling problem with blocking constraints and  multi-task flexibility is addressed. The authors prove a strong NP-hardness result, propose two mathematical models and  variable neighborhood search based heuristic algorithms.

Regarding the general case with any number of machines, in \cite{Ruiz2011} the authors investigate a flexible permutation flow shop, where operations must follow a predetermined order of assignment to machines and propose heuristic solution algorithms.
In~\cite{liao}, a permutation flow-shop in which one or more processors, consisting of one or more workers and/or facilities, is considered. Here, certain flexible processors execute their own set of tasks while also potentially assisting  other processors in performing their operations.

A different type of flexibility within flow shop and open shop models is considered in \cite{Knust2019}, where the authors examine \emph{pliable} jobs,  which are characterized by having predetermined total processing times. 
However, the specific processing times of operations constituting these jobs are not fixed and must be determined. 
They show that several variations of these problems featuring pliable jobs seem to be computationally simpler than their traditional counterparts.

A related concept known as \emph{resource flexibility} is explored in~\cite{danielsmazzola,daniels}, where the processing time of each operation is inversely proportional to the (discrete) amount of allocated resources. 
Our problem could fit within this framework by considering additional constraints on how resources are allocated to the machines.

\section{Problem statement and complexity}\label{sect:statement}

In this section we formally define the scheduling problem described in the introduction.
It will be called 
``Flow shop with Inter-stage fleXibility and Blocking'' ({\p}) and consists of the following elements:

\begin{itemize}
 \item A set $\mathcal{M}=\{M_1, M_2, \ldots,  M_m \}$ of $m$ machines  in a flow line, so that machine $M_k$ precedes machine $M_{k'}$ if $k<k'$.
 Buffers between machines are not available, i.e., intermediate storage capacity is considered zero.
 Thus, the so-called {\em  blocking} occurs, i.e.\ a job finished on machine $M_k$ blocks it until the next-stage machine $M_{k+1}$ is available for processing that job. 
 
 \item A set $\mathcal{J}=\{J_1, J_2, \ldots,  J_n \}$ of $n$ jobs, where each job $J_j$  consists of an ordered sequence $\langle o_{1j}, o_{2j}, \ldots, o_{q_j j} \rangle$ of $q_j \ge m$ operations to be processed in the given order. 
 Here, $o_{ij}$ denotes the $i$-th operation of job $J_j$. 
 \item  The $i$-th operation $o_{ij}$ of job $J_j$ 
 must be executed by exactly one machine chosen from a set $M(o_{ij})$.
 This set consists of either a single machine, i.e., $M(o_{ij}) =\{M_k\}$,
 or,  in case of a shiftable operation, of a pair of consecutive machines, i.e., $M(o_{ij})=\{ M_k, M_{k+1}\}$.
 The processing time of $o_{ij}$ on a machine $M_k\in M(o_{ij})$ is denoted by $p_{ij}^k$.
We assume that all the processing times are strictly positive.
 \end{itemize}

From the ordering of the operations and the flow line environment it follows that 
if $M(o_{ij})=\{M_k, M_{k+1}\}$ 
and $o_{ij}$ is processed by machine $M_{k+1}$, 
then none of the subsequent operations $o_{(i+1) j},o_{(i+2)j},\ldots,o_{q_j j}$ can be processed on the preceding machines $M_1, M_2,\ldots, M_k$.
In the remainder of the paper,
with no loss of generality, we assume that for every job $J_j$ and each machine $M_k$ there is a unique operation with $M(o_{ij})=\{M_k\}$. In fact, 
if there were several operations of this type, their processing times could be summed up in $p_{ij}^k$. In this case,
for the immediate preceding operation $o_{(i-1)j}$ either $M(o_{(i-1)j})=\{M_{k-1}\}$ or $M(o_{(i-1)j})=\{M_{k-1}, M_k\}$.
Similarly, for the succeeding operation $o_{(i+1)j}$ either $M(o_{(i+1)j})=\{M_{k+1}\}$ or $M(o_{(i+1)j})=\{M_{k}, M_{k+1}\}$. 
Note that, 
$(i)$ in this setting, as all operations processing times are strictly positive, every job ``visits'' all the $m$ machines of the flow line, in order (i.e., no machine can be skipped by any job); 
$(ii)$ the special case in which all operations have a single performing machine (i.e., $|M(o_{ij})| = 1$ for each $j\in\mathcal{J}$, $i=1,\ldots,q_j = m$) is the ``standard'' flow shop problem with blocking.

\begin{figure}
   \centering
       \includegraphics[width=\textwidth]{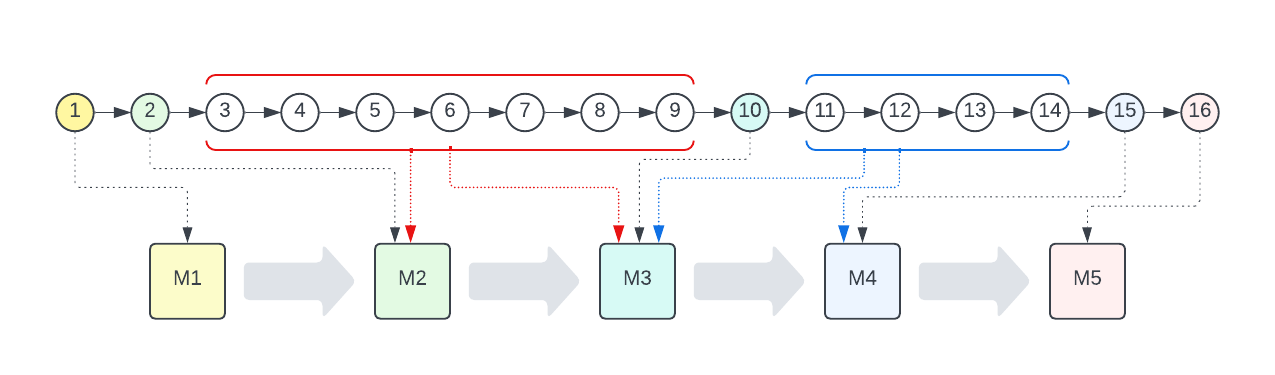}
      \caption{Structure of the operation-machine relationship in the real-world problem.}\label{fig:flowline}
\end{figure}

The structure of the problem requires an assignment decision for each operation of a job $J_j$ for which two machines are available.
A feasible assignment of the operations of $J_j$ to machines implies a partition of the set $\langle o_{1j}, \ldots, o_{q_j j}\rangle$ into $m$ subsets such that each subset consists of a certain number of consecutive operations 
assigned to a single machine able to process them.  
Figure~\ref{fig:flowline} illustrates all the above issues for an example with one job $J_j$ in which the sets $M(o_{ij})$, $i=1,\ldots,16$, are highlighted with different colors 
(for instance, operation $M(o_{10\, j}) = \{M_3\}$, while, for $i=11,\ldots,14$,  $M_{ij}=\{M_3,M_4\}$).

In the remainder of this paper, a feasible assignment is called \emph{assignment mode}.  It is straightforward to enumerate the set of all possible assignment modes for a job $J_j$, which we denote by $\mathcal{A}_j$.
If we call $n_k = |\{o_{ij} : M(o_{ij}) = \{M_k,M_{k+1}\}\}|$ the number of shiftable operations that can be processed either by machine $M_k$ or $M_{k+1}$,
then the number of feasible assignment modes for job $J_j$ is 
$$|\mathcal{A}_j| = \prod_{k=1}^{m-1} (n_k+1), \mbox{ which is } O\left(q_j-1 \choose m-1 \right).$$
For a given assignment mode $l\in\mathcal{A}_j$, 
we use $p_{j}^k(l)$ to indicate the (fixed) total processing time of all the operations of job $J_j$ assigned to machine $M_k\in\mathcal{M}$.  

As  already observed above, the flow line setting with blocking constraints implies that once we decide a processing order (i.e., a permutation) for the $n$ jobs, this remains the same on each machine, i.e., we are considering a variant of the {permutation flow shop}.
 
In conclusion, our special flow shop scheduling problem can be formulated as follows:
\begin{quote}
 {\sc Flow shop with Inter-stage Flexibility and Blocking} (\p): Given the set of linearly ordered machines $\mathcal{M}$, the job set $\mathcal{J}$, the processing times and a set of performing machines $M(o_{ij})$ for each operation $o_{ij}$ of every job $J_j\in\mathcal{J}$, we want to 
find 
$(i)$ a (feasible) assignment mode for each operation and
$(ii)$ a sequencing of the jobs,
such that blocking constraints hold and the makespan of the resulting schedule is minimized.  
\end{quote}

Based on our application and for ease of presentation, we assume that all jobs consist of the same number of operations with the same machine requirements. 
More precisely, for each pair of jobs $j, \ell\in\mathcal{J}$, $q_j=q_{\ell}=q$ and the $i$-th operation of jobs $j$ and $\ell$ can be executed on the same machines, i.e., $M(o_{ij})=M(o_{i \ell})$, $i=1,\ldots, q$. 
Note, however, that the models and algorithms presented in this paper do not depend on this assumption.

The complexity results (illustrated in Section~\ref{sec:literature}) for flow shop cases relevant to our study are summarized in Table~\ref{tab:complexity}. 
We are looking at problems with or without blocking constraints and flexibile operation-to-machine assignment. 
\begin{table}[h!]
    \centering
    \begin{tabular}{cccccc}
    \toprule 
         $n$ & $m$ & Fixed job sequence  & Blocking  & Flexibility  & Complexity\\
         \midrule 
         $n$ & 2  & yes & no  & yes  & NP-hard~\cite{Lin2017}\\
         $n$ & 2  & yes & yes & yes & Poly. Sect.\ref{sec:fixedsorting}  \\
         $n$ & 3  & yes & yes & yes & Open \\
         2   & $m$& no$^*$  & yes & yes  & Poly. Sect.\ref{sec:twojobs} \\   
         $n$ & 2  & no   & no  & no   & Poly.~\cite{Johnson54}\\
         $n$ & 2  & no   & yes & no   & Poly.~\cite{gg1964,Reddi1972} \\
         $n$ & Fixed $\ge 3$  & no & yes  & no  & NP-hard~\cite{Hall1996,papa80}\\
         $n$ & 2  & no   & yes & yes   & Open \\
    \bottomrule
    \end{tabular}
    \caption{Complexity cases ($^*$ With $n=2$ there are only two possible sequences.)}
    \label{tab:complexity}
\end{table}
While blocking simplifies the problem, flexibility introduces additional degrees of freedom and the problem becomes more complex. 
Its complexity remains open for the special case with $m=2$, flexibility, and blocking constraints.

\section{Special cases of \p}\label{sec:special_cases}

Here we consider two special cases of \p, namely the case in which there are only two jobs and the case in which there are only two machines and the sorting (i.e. the sequencing) of the jobs is fixed. We show that both problems can be solved in polynomial time.

\subsection{\p\ with two jobs}\label{sec:twojobs}

If only two jobs are to be processed in \p, the sequencing decision disappears as only two possible sequences have to be evaluated. Moreover, machine blocking becomes not relevant, as the minimum makespan takes the same value whether blocking holds or not. Still, the  operation-to-machine assignment problem remains non trivial. 

\begin{figure}[ht]
    \centering
    \includegraphics[width=0.7\linewidth]{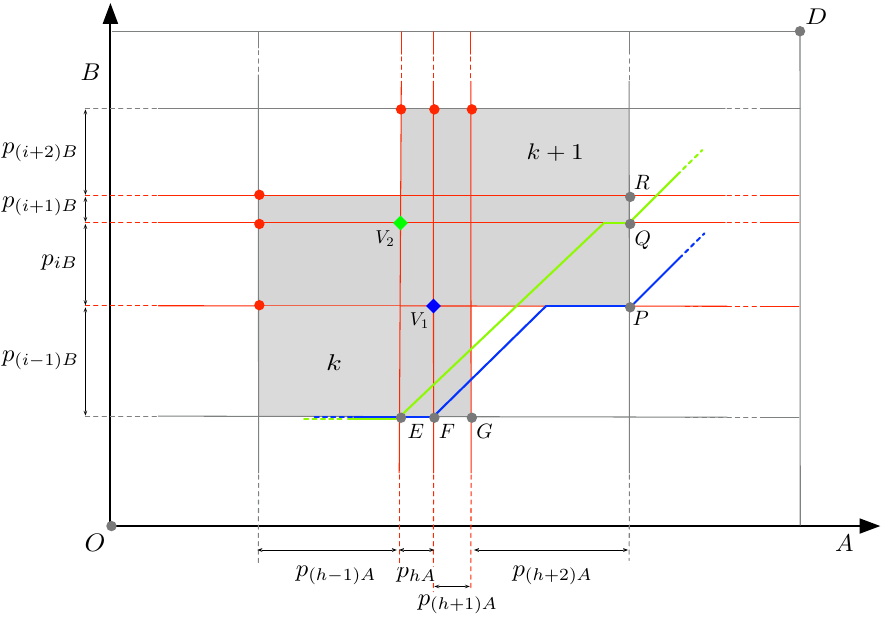}
    \caption{Grid representation of the $n=2$ case. 
    Here, the duration of a shiftable operation is the same on the two machines $M_k$ and $M_{k+1}$ ($p^{k}_{r\ell} = p^{k+1}_{r\ell} = p_{r\ell}$ for $r \ell =hA, iB, (h+1)A, (i+1)B$). Two paths associated to a sequence $A \prec B$ and different operation assignments ($F$ to $P$ and $E$ to $Q$), are depicted.
}
    \label{fig:grid}
\end{figure}

Hereafter, we show how a classical and well-known technique~\cite{AF55}, devised for the job shop scheduling problem with two jobs and based on a grid representation in a two-dimensional plane, can be extended to obtain an optimal assignment in polynomial time even with an arbitrary number of machines. 
For ease of presentation, the two jobs are denoted as $A$ and $B$,
where the horizontal (vertical) axis of the grid corresponds to job $A$ (job $B$).
We represent consecutive operations of each job by consecutive segments on their respective axes, 
where the lengths of the segments are equal to the duration of the operations. 
Starting from these segments we build a grid on the plane, as illustrated in~Figure~\ref{fig:grid}, with origin in the bottom-left point $O$, and an upper-right end point $D$.
If two operations are to be performed by the same machine $M_k$ and hence cannot be processed simultaneously, then we can associate to this pair an  ``obstacle'' on the grid. 
For instance, in Figure~\ref{fig:grid}, both operations $o_{(h-1)A}$ and $o_{(i-1)B}$ must be performed by machine $M_k$, so there is an obstacle---highlighted as a gray-colored rectangle---in correspondence to the segments $(h-1)$ and $(i-1)$ on the $A$ and $B$ axes, respectively.

Any allocation of the two jobs to the machines corresponds to an $(O,D)$-path on the grid consisting of horizontal segments, indicating that only job $A$ is processed, while job $B$ has to wait, vertical segments (the opposite), and 45-degree slanted segments, meaning that both jobs are being processed simultaneously. 
To preserve feasibility, such a path has to get around each encountered obstacle above or below, i.e., passing through top-left (NW) or bottom-right (SE) vertices of the obstacle. 
Intuitively, in order to save time, the path should take a diagonal segment unless an obstacle is met and, in that case the decision whether to get around the obstacle above  or below (through the NW or the SE corner) has to be taken.

Note that, since in {\p} only permutation schedules are feasible and there are only two job-sequences, then the paths only go through the SE corners (for the sequence $A \prec B$), or only through the NW corners (sequence $B \prec A$) of the obstacles. In the remainder of this Section, we are only considering the sequence  $A \prec B$, as the case $B \prec A$ is clearly symmetric.

Differently from the model proposed in~\cite{AF55}, for {\p} the size of the obstacles, and hence the actual positions of NW and SE corners, depends on the specific operation-to-machine assignments. Figure~\ref{fig:grid} illustrates this concept.
Here, operations $o_{hA}$, $o_{(h + 1)A}$, $o_{iB}$ and $o_{(i+1)B}$  are shiftable, i.e., can be processed by both machines $M_k$ and $M_{k+1}$, while operations $o_{(h-1)A}$ and $o_{(i-1)B}$ can be processed only by machine $M_k$ and $o_{(h+2)A}$ and $o_{(h+2)B}$ only by $M_{k+1}$. 
In the picture, every feasible operation assignment is uniquely identified by a pair of points $(s,r)$, with $s\in\{E,F,G\}$ and $r\in\{P,Q,R\}$ determining the assignments of the operations of, respectively, job $A$ and $B$. 
This way, points $s$ and $r$ also correspond to alternative positions of the SE corners of obstacles associated to machine $M_k$ and $M_{k+1}$ respectively, which in turn would correspond to different paths between $s$ and $r$.
For instance, if operation $o_{hA}$ is assigned to machine $M_k$ while $o_{(h+1)A}$, $o_{iB}$, and $o_{(i+1)B}$ to $M_{k+1}$, then $s=F$ and $r=P$ and we would choose the dark colored path that goes through these two points. 
In this case the upper-right corner and bottom-left corners of the two obstacles produced by machine $M_k$ and $M_{k+1}$ would coincide with point $V_1$.
Instead, the light colored path passing through points $E$ and $Q$, corresponds to assigning operation $o_{iB}$ to machine $M_k$, and operations $o_{hA}$, $o_{(h+1)A}$, and $o_{(i+1)B}$ to machine $M_{k+1}$ and the two obstacles meet at point $V_2$.

It can be easily shown that an $(O,D)$-path representing a minimum makespan schedule of \p\ corresponds to a shortest path in a suitable acyclic graph 
in which nodes correspond one-to-one with all SE corners of the obstacles representing any possible assignment (e.g., points $E$, $F$, $G$, $P$, $Q$, $R$ in Figure~\ref{fig:grid}) plus points $O$ and $D$.
In the graph, 
an arc exists between nodes $u$ and $v$ if a 45-degree line starting from the SE corner corresponding to $u$ on the grid and extending towards the top right intersects the obstacle whose SE corner corresponds to $v$ (or, if $v = D$, it reaches the border of the grid.) The length of arc $(u,v)$ equals the minimum time required to perform the operations from node $u$ to node $v$, which is the maximum value between the lengths, along the $x$ and $y$ axes, between the two nodes. 
For instance, the length of arc (i.e.\ the part of the dark path) between $F$ and $P$ in~Figure~\ref{fig:grid} is given by $\max\{p^{k+1}_{(h+1)A}+p^{k+1}_{(h+2)A}, p^k_{(i-1)B}\}$. 

Note that, in the figure the processing time of a shiftable operation is independent of the machine performing it. The general case, where this restriction does not apply, can be easily handled in the graph model by appropriately adjusting the arc lengths.

As observed in~\cite{Brucker88}, building the graph 
and computing a shortest path (i.e., the schedule makespan) can be done in $O(t \log t)$ time, where $t$ is the number of incompatible pairs.
In conclusion, the above algorithm for \p\ with $n=2$ jobs runs in polynomial time (using, e.g., $A^*$ algorithm~\cite{Agnetis1995,Nicosia2003}), even if the number of machines and number of shiftable operations is part of the input. 

\subsection{\p\ with two machines and fixed sequencing of the jobs}\label{sec:fixedsorting}
Even if the job sequence is given and fixed, minimizing makespan in a flow shop with only two machines and with only three operations per job (and one should decide which machine the intermediate operation has to be assigned to) is binary NP-hard~\cite{Lin2017}. 
Thus inter-stage flexibility increases the difficulty of the problem from a complexity point of view. 
However, if we consider blocking constraints in addition, the problem actually becomes easier.
We will show this result, which is somewhat counterintuitive on first sight, in the remainder of this section.

We consider the case with $m=2$ machines and a fixed sequence of the jobs, numbered from $J_1$ to $J_n$.
In this case, each job $J_j$ with $q_j$ operations 
must be processed in one of $|\mathcal{A}_j|$ assignment modes, where an assignment mode $l_j \in \mathcal{A}_j$ corresponds to the possibility of processing the first $l_j$ operations on machine $M_1$ and the remaining $q_j-l_j$ operations on $M_2$. Since, for any assignment mode, operation $o_{j1}$ is processed on $M_1$ and operation $o_{j q_j}$ on $M_2$, there is  
$l_j=1,\ldots, q_j-1$.

As pointed out in~\cite{Hall1996}, in general, flow shop with blocking constraints is a relaxation of the corresponding problem with no-wait constraints. However, it is easy to see that when there are only two machines the values of optimal solutions of the two versions of the problem are equal.
Hence, w.l.o.g.\ for every optimal solution of a flow shop problem with two machines and blocking constraints, we can assume that 
every job starts as late as possible on $M_1$, so that after its completion on $M_1$ it can be immediately moved to $M_2$. 

\begin{figure}[b]
    \centering
    \includegraphics[width=0.7\linewidth]{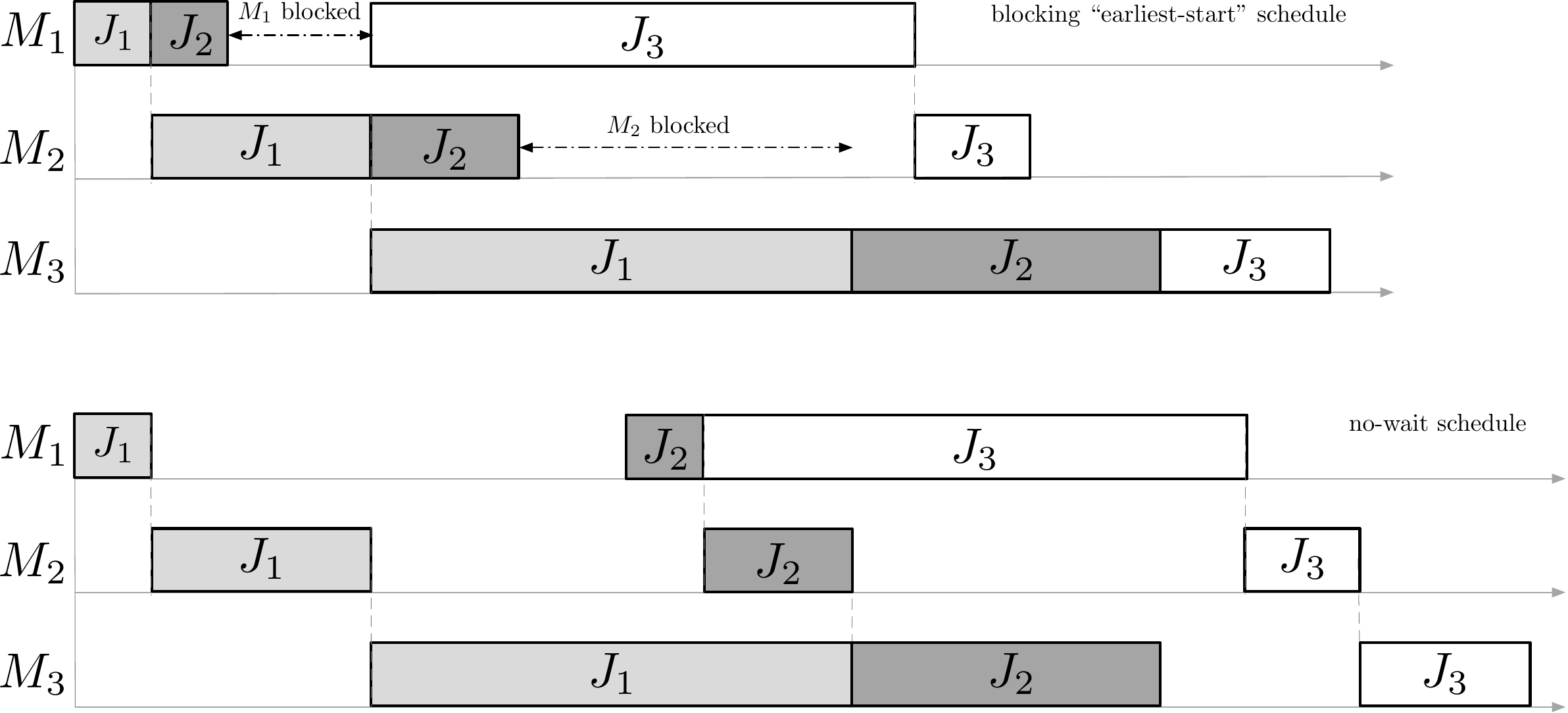}
    \caption{Minimum-makespan schedules: flow shop with blocking and the corresponding no-wait version}
    \label{fig:waiting}
\end{figure}

\begin{lemma}\label{th:latestart}
When $m=2$ (i.e., there are only two machines), there always exists an optimal solution of {\p} such that, for every job, the completion time on $M_1$ is equal to the starting time on $M_2$.
\end{lemma}
\begin{proof}
If some job $J_j$ finishes its last operation on $M_1$ before starting its first operation on $M_2$, we can always postpone the starting time of $J_j$ on $M_1$ such that the condition of the Lemma is fulfilled. 
This has no consequences for job $J_{j+1}$ since $J_j$ blocks $M_1$ until the start of its processing on $M_2$.
\end{proof}
It should be noted that Lemma~\ref{th:latestart} does not hold anymore with $m\ge 3$ machines as illustrated by the example in Figure~\ref{fig:waiting}.
There, it is beneficial to start $J_2$ early on $M_2$ and thus permit $J_3$ to start earlier on $M_1$.
This situation cannot occur for two machines.

Based on the structural property of Lemma~\ref{th:latestart} we can establish the computation of the objective function value of \p.
Assume that for every job $J_j$ a certain assignment mode $l_j \in \mathcal{A}_j $ is chosen yielding total processing times of
$\pp_j^1 := \sum_{i=1}^{l_j} p_{i j}^1$ on $M_1$ and
$\pp_j^2 := \sum_{i=l_j+1}^{q_j} p_{i j}^2$ on $M_2$.

If $\pp_j^1 \geq \pp_{j-1}^2$ then job $J_{j-1}$ is already finished on $M_2$ when $J_j$ moves from $M_1$ to $M_2$ leaving  idle time of $\pp_j^1 - \pp_{j-1}^2$ on $M_2$.
If $\pp_j^1 < \pp_{j-1}^2$ then there is necessarily idle time of $\pp_{j-1}^2 - \pp_j^1$ on $M_1$ while $J_{j-1}$ is being processed on $M_2$.
In this case, the contribution of $J_j$ to the total completion time consists of the processing time on $M_1$ plus the idle time on $M_1$ which sums up to $\pp_{j-1}^2$.
Therefore, the overall completion time can be written as
\begin{equation}\label{eq:m2compltime}
   C_{\max} = \pp_1^1 + 
   \sum_{j=2}^{n} {\max\{\pp_j^1,\, \pp_{j-1}^2\}}
   +\pp_n^2\,.
\end{equation}

\begin{figure}[t!]
    \centering
    \includegraphics[width=0.8\linewidth]{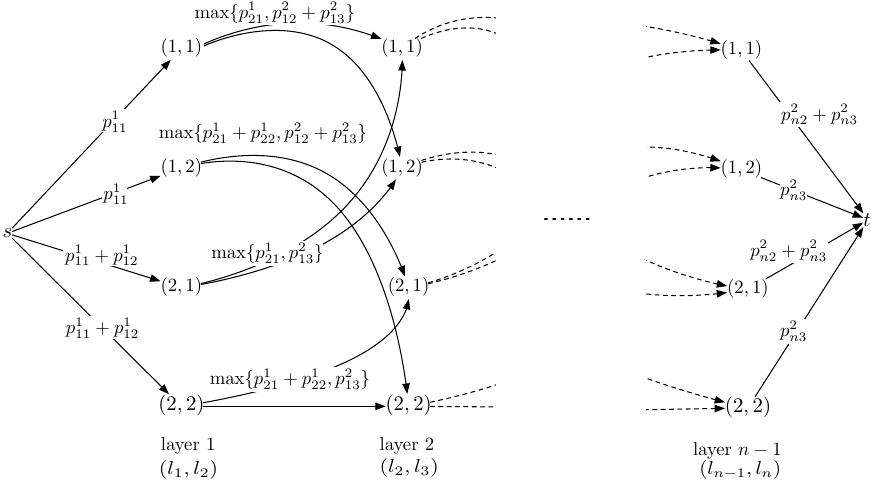}
    \caption{Digraph for makespan minimization}\label{fig:stategraph}
\end{figure}

To determine the optimal assignment mode for each job and thus minimize the completion time $C_{\max}$ for \p\ we introduce the following directed graph model $G=(N,A)$.
The node set $N$ consists of $n-1$ layers, one layer for each job $J_1,\ldots, J_{n-1}$.
In each layer $j$ there are $(q_j-1) \cdot (q_{j+1}-1)$ nodes.
Each such node represents the decision to choose a certain pair of assignment modes $(l_j, l_{j+1}) \in \mathcal{A}_j \times \mathcal{A}_{j+1}$ for consecutive jobs $J_j$ and $J_{j+1}$.
In addition, there is a source node $s$ and a sink node $t$.
Arcs in $A$ connect each node in layer $j$ to all those nodes in layer $j+1$ ($j=1,\ldots, n-1$) where the assignment mode $l_{j+1}$ chosen for job $J_{j+1}$ in the two nodes coincide.
Furthermore, there are arcs from $s$ to every node in layer $1$, and from every node in layer $n-1$ to node $t$.

All 
arcs from node $(l_{j-1}, l_{j})$ of layer $j-1$ to nodes  $(l_{j}, l_r)$ in layer $j$, with $j=2,\ldots, n-1$, $r=1, \ldots, q_{j+1} - 1$, are given a weight value representing the increase of the overall completion time $C_{\max}$ implied by the assignment modes $l_{j-1}$ and $l_{j}$ chosen for jobs $J_{j-1}$ and $J_{j}$, namely $\max\{\pp_j^1,\, \pp_{j-1}^2\}$  according to Equation~\eqref{eq:m2compltime}.

For the arcs from $s$ entering  a node in layer $1$ referring to job $J_1$ we have weights $\pp_1^1$, depending on the assignment mode represented by that node.
Finally, the arcs from layer $n-1$ into $t$ have a weight 
$\max\{\pp_n^1,\, \pp_{n-1}^2\}+\pp_n^2$ indicating the contribution of the final job $J_n$ under the respective assignment mode chosen for $J_n$ in the node in layer $n-1$.

Figure \ref{fig:stategraph} illustrates the above described layered graph for the special case in which $q_j = 3$ for all jobs $J_j\in\mathcal{J}$ and only the intermediate operation is shiftable, i.e., it can be assigned alternatively to the first ($l_j=2$) or the second ($l_j=1$) machine. 
In this case, for the $j$-th job $J_j$ we have: 
$\pp_j^1 = p_{j 1}^1 + (l - 1)p_{j 2}^1$ and
$\pp_j^2 = (2 - l)p_{j 2}^2 + p_{j 3}^2$.

A shortest $(s,t)$-path in  graph $G$ corresponds to a decision with minimum makespan $C_{\max}$.
Thus, we conclude:
\begin{theorem}
When $m=2$ (i.e., there are only two machines) and the sequence of jobs is fixed, \p\ can be solved in polynomial time for arbitrary many operations.
\end{theorem}
\begin{proof}
The number of nodes in the graph $G$ can be bounded by $O(n\, q_{\max}^2)$, where $q_{\max}:=\max_j q_j$.
This is polynomial in the input size.
\end{proof}

\section{Heuristic use of MIP models}\label{sec:heur_ilp}
In this section, we introduce a series of matheuristic algorithms based on two of the mixed integer programming (MIP) models for \p\ we presented in \cite{us_ods2023}. 
Matheuristics, which combine mathematical programming techniques with heuristic approaches to solve optimization problems, have gained significant popularity over the last decade due to the significant advancements in computational capabilities. 
In the field of scheduling, matheuristics have been extensively studied, particularly for tackling computationally challenging shop scheduling problems, see, for instance, \cite{Gao2022,LIN2016} for no-wait flow shop,  
\cite{Yang2021} for flow shop problems with time windows, and \cite{Pastore2022} for open shops with time dependent setup times.

\begin{figure}[htp]
    \centering
    \includegraphics[width=1.0\linewidth]{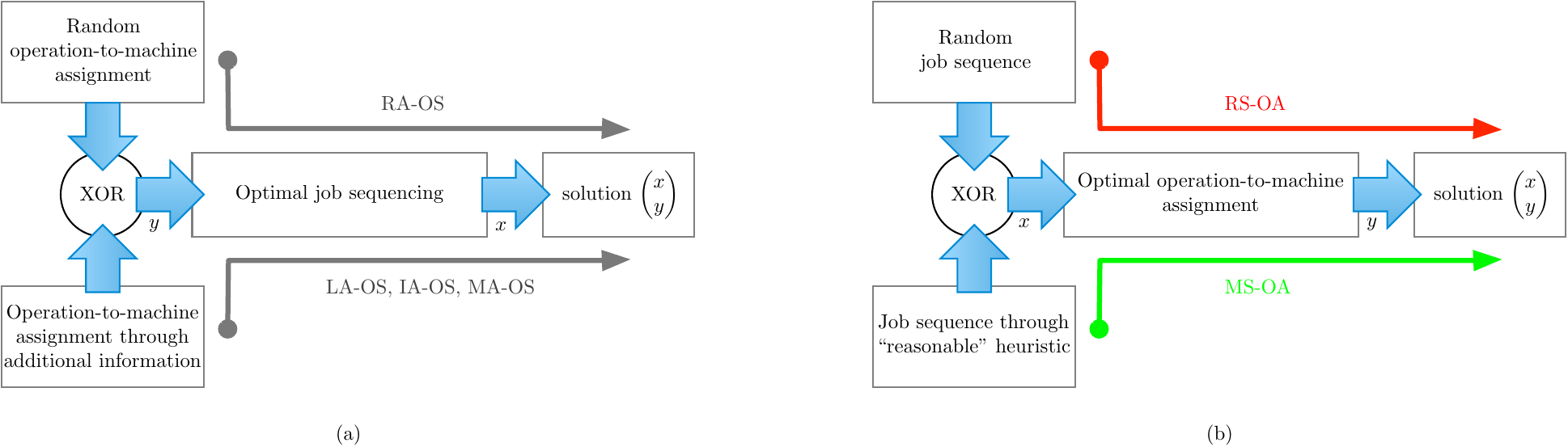}
    \caption{Scheme of heuristic approaches based on MIP1: fix assignment first (a) and fix sequence first (b)}
    \label{fig:HeurScheme}
\end{figure}

The matheuristics that we propose follow  distinct approaches based on the order in which the two main decision tasks—operations assignment modes and job sequencing—are tackled. 
In fact, one of the MIP that we consider is well suited for separating the two decisions and is therefore used in Sections \ref{sec:asgnfirst} and \ref{sec:seqfirst} to formulate different algorithmic strategies, which are illustrated in Figure~\ref{fig:HeurScheme}.
In Section \ref{sec:holistic}, we present an alternative method that uses a different MIP model and addresses both decision tasks simultaneously through an iterative rounding procedure.

\subsection{MIP models}\label{sec:mip}

In \cite{us_ods2023}, the authors propose four different mixed integer linear programs for \p. 
In two of them so-called ``positional''  variables are used to define the sequencing of the jobs.
The assignment of operations to machines is addressed in two different ways.
In the \emph{explicit assignment model} a dedicated set of variables is used to \emph{explicitly} model the operation-to-machine assignments. 
In the \emph{implicit assignment model} all possible assignment modes $l\in\mathcal{A}$ are computed in a preprocessing phase together with the resulting job processing times $p_{j}^k(l)$ on machine $M_k\in\mathcal{M}$, $J_j\in\mathcal{J}$, $l\in\mathcal{A}$. 
Then the model considers how to sequence the jobs in the most effective assignment mode.
Hereafter, we present only the two ``positional'' models of \cite{us_ods2023}, as they were shown to be the best-performing ones. They will be used in the matheuristics discussed below.

The {\em positional variable model with explicit assignment}, from now on denoted as {\bf MIP1}, considers binary variables $ x_{jh} $ indicating whether job $J_j$ is processed in the $h$-th position and $y_{ijk}$ indicating if operation $o_{ij}$ of job $J_j$ is assigned to machine $M_k$.
$s_{hk}$ are continuous variables indicating the starting time of the job in position $h$ on machine $M_k$;
$P_{hk}$ are continuous variables indicating the total processing time of the operations which are assigned to machine $k$ and belong to the job in position $h$.
Moreover, we denote by $\mathcal{H}=\{1,2 \ldots, n\}$ the set of all positions.
\begin{align}
  \min &\  C_{\max}= s_{nm} + P_{nm} \ &   \label{cmax_01} \\
   s.t.\     &\ \sum\nolimits_{h\in\mathcal{H}} x_{jh} = 1 & J_j \in \mathcal{J}\label{assign1} \\ 
         &\ \sum\nolimits_{J_j\in J} x_{jh} = 1 & h \in \mathcal{H} \label{assign2}\\ 
         &\ \sum\nolimits_{M_k \in M(o_{ij})} y_{ijk} = 1 &  i=1,\ldots, q_j, J_j  \in \mathcal{J} \label{assign3}\\ 
         & y_{i+1, j k} + y_{ijk'} \leq 1   
         \begin{split}  & i=1,\ldots, q_j, J_j \in \mathcal{J},  \\ 
         & M_k \in M(o_{i+1,j}),  M_{k'} \in M(o_{ij})\: :\: k < k' 
         \end{split}\label{cutdefinitionbis}\\
           &\ s_{h, k+1} \geq s_{hk} + P_{hk} &  M_k \in \mathcal{M}\setminus \{ m \} , h\in\mathcal{H}   \label{precmachine} \\
          &\  s_{h+1, k} \geq s_{hk} +  P_{hk} & M_k \in \mathcal{M} , h\in\mathcal{H}\setminus\{n\}   \label{precjob}\\                           
            &\ s_{hk} \geq s_{h-1, k+1}  &  M_k \in \mathcal{M} \setminus \{ m \} , h \in\mathcal{H}\setminus\{1\}   \label{blocking}  \\
              \begin{split}& P_{hk} \geq \sum\nolimits_{i : M_k\in M(o_{ij})} p_{ij}^k y_{ijk} \\ 
             & \quad \quad \quad \quad \quad \ - B_{jk}(1-x_{jh})
          \end{split}
             &   h \in\mathcal{H}, M_k \in \mathcal{M}  , J_j \in \mathcal{J}  \label{processingonmachinek} \\ 
             &\  x_{jh}, y_{ijk}  \in \{0,1\} & i=1,\ldots, q_j, J_j \in \mathcal{J}, h\in\mathcal{H}, M_k \in \mathcal{M} \label{binary}\\
            &\ s_{hk}, P_{hk} \ge 0  & M_k \in \mathcal{M} , h\in\mathcal{H}   \label{cont}
            \end{align}
The ``large'' constant $B_{jk}$ in \eqref{processingonmachinek} may be set equal to $\sum_{i : M_k\in M(o_{ij})} p_{ij
}^k$.

The makespan is represented by \eqref{cmax_01}, that is the completion time of the last job on the last machine.
Constraints~\eqref{assign1} and \eqref{assign2} assign exactly one job to each position and \eqref{assign3} assigns each operation to exactly one machine.
Respecting the sequence of operations on the machines is enforced by \eqref{cutdefinitionbis}.
The correct sequence of jobs w.r.t.\ machines is implied by constraints~\eqref{precmachine} and \eqref{precjob}.
Expression \eqref{blocking} imposes the blocking constraints, enforcing that the job in position $h$ can only start being processed on machine $M_k$ when the preceding job in position $h-1$ has moved on and started its processing on machine $M_{k+1}$.
The definition of the processing time of the job placed in position $h$ on machine $M_k$ is given by \eqref{processingonmachinek}.

\medskip
The {\em positional variable model with implicit assignment}, from now on denoted as {\bf MIP2}, turned out in \cite{us_ods2023} to be the best performing one among all four models.
It uses the following sets of variables: 
binary variables $ x_{jhl} $ indicate whether job $J_j$ is processed in the $h$-th position according to operation assignment $l$;
as before continuous variables $s_{hk}$ indicate the starting time of job in position $h$ on machine $M_k$, and $C_{\max}$ is the schedule makespan in \eqref{ob:cmax}.
\begin{align}
  \min &\ C_{\max} = s_{n m} + \sum\nolimits_{\ncoppa{J_j\in\mathcal{J}}{l \in \mathcal{A}}} p_{j}^m(l) x_{jnl}  \ & \label{ob:cmax}   \\
   s.t.\       &\ \sum\nolimits_{\ncoppa{h\in\mathcal{H}}{l \in \mathcal{A}}} x_{jhl} = 1 & J_j \in \mathcal{J}\label{LPv1} \\ 
             &\ \sum\nolimits_{\ncoppa{J_j\in \mathcal{J}}{l \in \mathcal{A}}} x_{jhl} = 1 & h \in \mathcal{H}  \label{LPv2}\\ 
           &\ s_{h, k+1} \geq s_{h k} + \sum\nolimits_{\ncoppa{J_j\in \mathcal{J}}{l \in \mathcal{A}}} p_{j}^k(l)\, x_{jhl} &  M_k \in \mathcal{M} , h \in \mathcal{H}   \label{LPv4} \\
          &\  s_{h+1, k} \geq s_{h k} + \sum\nolimits_{\ncoppa{J_j\in \mathcal{J}}{l \in \mathcal{A}}} p_{j}^k(l)\, x_{jhl} &  M_k \in \mathcal{M} , h \in \mathcal{H}  \label{LPv5}\\                      &\ s_{h k} \geq s_{h-1, k+1}  &  M_k \in \mathcal{M} , h \in \mathcal{H} \setminus \{1\}   \label{LPv6}  \\
            &\ x_{jhl} \in \{0,1\} &   J_j \in \mathcal{J}, h\in \mathcal{H},  l \in \mathcal{A} \label{LPv8}\\
            &\ s_{h k}
            \ge 0  & M_k \in \mathcal{M} , h \in \mathcal{H}    \label{LPv10}
            \end{align}
Constraints \eqref{LPv1} assign exactly one position and one assignment mode to each job, while constraints \eqref{LPv2} assign exactly one job  and one assignment mode to each position.
Constraints \eqref{LPv4} ensure that each job in any position $h$ cannot start its processing on a machine $M_{k+1}$ before it has finished its operations on machine $M_k$. 
Constraints \eqref{LPv5} ensure that on any machine $M_k$ the job in position $h+1$ can only start its processing after the operations of the job in position $h$ were completed on machine $M_k$. 
Constraints \eqref{LPv6} are  the blocking constraints, enforcing that the job in position $h$ cannot  start being processed on machine $M_k$ unless the preceding job in position $h-1$ has moved on and started its processing on machine $M_{k+1}$.

\subsection{Matheuristics fixing assignment first}\label{sec:asgnfirst}
Here, we describe MIP-based heuristic algorithms which determine first an assignment of operations---thus fixing the $y_{ijk}$ variables---and then compute the \emph{optimal} sequence represented by variables $x_{ij}$ by solving a restricted instance of MIP1. 
See the sketch in Figure~\ref{fig:HeurScheme}(a) for a graphical representation of the solution approach.

We have tested four different methods to decide the initial operation assignment which are listed hereafter.
\begin{description}
    \item[RA-OS] A first straightforward approach consists in randomly assigning operations to feasible machines. The resulting matheuristic is called ``{Random Assignment - Optimal Sequence}" (RA-OS). Note that, in doing so, all the precedence relations detailed in Section~\ref{sect:statement} are to be taken into account. 
    For instance, if job $J_j$ is comprised of a sequence of operations $\langle o_{1j},\ldots,o_{q_j j}\rangle$ and operation $o_{ij}$ can be executed by both machines $M_k$ and $M_{k+1}$, then fixing $y_{ijk}:=1$ implies $y_{\ell jk}:=1$ for all $1\le \ell < i$; if otherwise  $y_{ij, k+1}:=1$ then $y_{\ell j, k+1}:=1$ for all $1 < \ell \le q$.
    \item[LA-OS] To start with a more promising assignment we can exploit information from the LP-relaxation of MIP1: let $\bar{y}_{ijk}$ be an optimal value of the $y_{ijk}$ variable for such a relaxation. Assuming that operation $i$ of job $j$ can be executed either on machine $M_k$ or machine $M_{k+1}$, we assign it to machine $k$ if $\bar{y}_{ijk} \geq 0.5$, and to machine $M_{k+1}$ otherwise. This matheuristic is called ``{LP-rounding Assignment - Optimal Sequence}" (LA-OS).
    \item[IA-OS] Information coming from the LP-relaxation can be exploited even further by proceeding as follows.
    After solving the LP-relaxation of MIP1 we choose the largest value $\bar{y}_{ijk}$ and set it to $1$, which obviously forces a number of other $y$-values to $0$. 
    After fixing these values we iterate the procedure, by solving again the LP with the remaining free variables (keeping fixed the others). 
    This will be called ``{Iterative-LP assignment - Optimal Sequence}" (IA-OS).
    \item[MA-OS] An alternative way without using the LP-relaxation consists of determining, for each job $J_j$, the values of the assignment such that the total processing time of $J_j$ is minimum, i.e., the locally optimal assignment mode for $J_j$. 
    Given these values of assignment variables $y^{\min}_{ij}$, we compute again an optimal sequence of jobs by solving MIP1 to optimality yielding variable values $x^{\min}_{ij}$. 
    This approach is called ``{Min-time Assignment - Optimal Sequence}" (MA-OS).
\end{description}

\subsection{Matheuristics fixing sequence first}
\label{sec:seqfirst}

Exchanging the order in which the two decisions are taken, we get a different approach for a matheuristic. 
Thus, the idea of the solution approaches that we describe below is to first fix the job sequencing and then optimize the operation-to-machine-assignment. 
The job sequencing is defined by fixing the $x_{ij}$ variables in MIP1 and afterwards, the resulting  restricted MIP with only $y$ variables is solved to optimality. 
Also in this case, we propose different ways to determine the initial sequence (see the illustration in Figure~\ref{fig:HeurScheme} (b)) yielding the following algorithms.
\begin{description}
    \item[RS-OA] An easy way to tackle the first decision task, is generating a random sequence of jobs and fixing the resulting values of the $x_{ij}$ variables accordingly. 
    Then we go back to model MIP1 and optimally solve the restricted model on the $y_{ijk}$ free variables with fixed $x_{ij}$, hence determining the optimal operation-to-machine assignment. 
    The resulting matheuristic is called ``{Random Sequence - Optimal Assignment}" (RS-OA). 
    \item[MS-OA] Clearly, one could expect to improve the outcome of the solution approach by choosing the starting sequence of jobs in a more considerate way. 
    However, the optimization of the job sequencing should be based on job processing times, which again depend on the assignment of operations. 
    To break this loop, we use the 
    sequence output by the matheuristic MA-OS described above (hence assuming the minimum total processing time for each job) and fix the corresponding optimal starting sequence accordingly, i.e., fixing $x_{ij}:=x^{\min}_{ij}$. 
    The decision on the assignment of operations is then revised, again by solving the restricted model on the $y_{ijk}$ free variables with fixed $x^{\min}_{ij}$ variables in the MIP1 model. 
    We call this two-step procedure ``{Min-time Sequence  - Optimal Assignment}'' (MS-OA). 
\end{description}

\subsection{Matheuristic with a holistic approach}\label{sec:holistic}

As an alternative to the two previous types of matheuristics with a hierarchical approach to the two decision tasks, we now introduce a more comprehensive strategy and aim at gaining advantages from considering the problem as a whole, trying to exploit information from the second model MIP2 as much as possible.
In particular, we derive information from the optimal solution of the LP-relaxation of 
MIP2 (positional model with implicit assignment), which combines sequencing and assignment in one set of variables $x_{jhl}$.

In our approach called ``{Sequential-fixing-with-threshold}'' ({\bf SFT}), we repeatedly solve the latest linear program (initially, the LP-relaxation of MIP2).
Let $\bar x$ be the optimal solution of the current LP-relaxation and sort the $x$ variables in non-increasing order of the component values $\bar x_{jhl}$. 
Given a threshold $\phi \in [0,1]$ and an integer parameter $1\le r\le n$, consider the $r$ variables with largest $\bar x_{jhl}$ values:
If $\bar x_{jhl} \geq \phi$, then the variable $x_{jhl}$ is permanently fixed to $1$ in the current and in all succeeding LPs.
The resulting linear program is again solved to optimality and the procedure is iterated. 
If all variables are below the threshold $\phi$, the remaining MIP model (i.e, MIP2 with the remaining free variables, keeping fixed the previous ones) is solved yielding an optimal solution for the restricted problem.

One may observe that when $r$ is small, SFT likely goes through a considerable number of iterations consuming rather high running times with respect to the case in which a larger number $r$ of variables can be simultaneously fixed at each iteration.
For instance, if $r=1$ at most one variable $x_{ijk}$ can be fixed in each LP iteration. Instead, if $r=n$, all variables whose values in $\bar x$ reach the threshold $\phi$ are simultaneously fixed in the same LP iteration: In this case, we can expect shorter computation times, but at the cost of less accurate values for the $x$ variables, which may ultimately reduce the quality of the final solution.   

Depending on the choice of the number of fixable variables $r$ and the threshold $\phi$, we obtain different variants for the SFT matheuristic 
 denoted as {\bf SFT$(r,\phi)$}. 
Differences in performance between variants of SFT will be described in Section~\ref{sec:comparison_heur}.

\section{Constructive Insertion Heuristic}\label{sec:constrheu}

We propose a constructive heuristic for \p\ based on the following principle:
We start with an empty schedule and iteratively add jobs to build a complete schedule.

Consider the jobs in a given order. 
Let $\sigma_j$ denote  the current partial schedule corresponding to the sequencing and operation assignment of the first $j$ jobs and let $\sigma_j$ be empty when $j=0$. The idea of the algorithm is, for each job $J_j$ (with $j=1, 2, \ldots , n$),  to consider its inclusion in the current partial sequence defined by $\sigma_{j-1}$ by trying all $j$ possible positions for placing $J_j$.

For each insertion position, 
we heuristically determine  the assignment of the operations for $J_j$.
At the end of each such step, we choose for job $J_j$ the insertion position 
and the operation assignments 
which minimize the overall makespan of the first $j$ jobs. 
This way a new partial schedule $\sigma_j$ is generated. 
When considering an insertion position $h$, job $J_j$ is inserted directly after the job currently at position $h$ in $\sigma_{j-1}$, with $h=0,1, \ldots, j-1$. 
After the insertion, $J_j$ is in position $h+1$ in the new schedule $\sigma_{j}$ and our heuristic leaves all the assignments of operations for jobs on positions $1, 2, \ldots, h$ unchanged, while the operation assignments for the remaining jobs, scheduled later than $J_j$, 
are computed as follows.
For each job \( J_\ell \) in positions \( h+1, \ldots, j \), and for each machine \( M_k \), from \( k = 1 \) to \( m-1 \), if there are shiftable operations that can be assigned to either machine \( M_k \) or \( M_{k+1} \) (i.e., \( n_k \geq 1 \)), we assign these operations in a way that minimizes the completion time of the operation \( o_{i \ell} \), where \( M(o_{i \ell}) = \{M_{k+1}\} \) (the first non-shiftable operation of \( J_\ell \) on machine \( M_{k+1} \)).
\begin{algorithm}[htbp]
\begin{algorithmic}[1]
\STATE {\bf Input}: job set $J_1,J_2,\ldots,J_n$
\STATE $\sigma \longleftarrow$ empty schedule  
\FOR{$j=1\ldots n$} 
       \FOR{$h=j-1$  down to $0$} 
            \STATE $\sigma_{jh} \longleftarrow$ \textsc{Insert $(\sigma,j,h)$} 
            \COMMENT{insert $J_j$ in $\sigma$ after the job currently at position $h$}
            \STATE  $C(h) \longleftarrow$ completion time of the last job of $\sigma_{jh}$ on $M_m$
     \ENDFOR
    \STATE $\sigma \longleftarrow \sigma_{jh^*}$ with $h^* = \arg\min_{h = 0,1,\ldots,j-1}\{C(h)\}$
\ENDFOR
\STATE {\bf Output}: $\sigma$ \COMMENT{Schedule of the $n$ jobs comprising the operation assignment and the sequencing of all jobs} 
\end{algorithmic}
     \caption{Insertion Heuristic}
    \label{alg:constr}
\end{algorithm}  

\begin{algorithm}[htbp]
\begin{algorithmic}[1]
\STATE {\bf Input}: Schedule $\sigma$, job index $j$, position $h$ 
\STATE Let $\langle J_{[1]}, J_{[2]} \ldots J_{[j-1]}\rangle$ be the sequencing of the jobs in $\sigma$ \COMMENT{$J_{[p]}$ indicates the job in position $p$ in the sequence}
\STATE Set the sequencing of $\sigma'$ to be $\langle J'_{[1]}, J'_{[2]}, \ldots, J'_{[j]}\rangle = \langle J_{[1]}, \ldots  J_{[h]},  J_j , J_{[h+1]}, \ldots , J_{[j-1]}\rangle $
\STATE Set the operation assignment of jobs in  $\sigma'$ at positions $1,\ldots, h$  the same as in $\sigma$ 
  \STATE{Compute the operation assignment for jobs at  positions $h+1 \ldots j$ in $\sigma'$: }   
  \FOR{$p=h+1 \ldots j$}
       \STATE Let $J_{\ell}= J'_{[p]}$ be the job at position $p$ in $\sigma'$
        \FOR{$k=1 \ldots m-1$}
            \STATE Let $O_\ell^k:=\{o_{i \ell} \mid M(o_{i \ell})=\{M_k, M_{k+1}\}\}$ \COMMENT{set of shiftable operations of job $J_\ell$ that can be executed by both machines $M_k$ and $M_{k+1}$}
            \STATE Assign operations $O_\ell^k$ to $M_k$ or $M_{k+1}$ so that the resulting completion time of the unique operation $o_{i \ell}$ of $J_\ell$ with $M(o_{i \ell})=\{M_{k+1}\}$ is minimized. 
    \ENDFOR 
    \ENDFOR 
\STATE {\bf Output}:  $\sigma'$ \COMMENT{Schedule in which $J_j$ is inserted in $\sigma$ after the job at position $h$}
\end{algorithmic}
     \caption{Function {\sc Insert}}
    \label{alg:constr2}
\end{algorithm}  

A sketch of the above described procedure is given in Algorithm~\ref{alg:constr} and~\ref{alg:constr2}.

The running time of the Insertion Heuristics can be bounded as follows.
Every execution of Function {\sc Insert} takes $O(nm)$ time.
It is called $n^2$ times from the main body of the heuristic.
These calls dominate the total running time which yields an overall running time complexity of $O(n^3 m)$.

\section{Computational Results}
\label{sec:results}

In this section we present the results of a computational campaign carried out to assess the
computational efficiency of the different solution approaches presented in Sections~\ref{sec:heur_ilp} and~\ref{sec:constrheu}.
The experiments have been run on an 11th Gen Intel® Core™ i7-1165G7 @ 2.80GHz × 8 with 16GB of available RAM, running Debian 11. 
The procedure was implemented in Python; the optimization problems were solved with Gurobi 11.0.1. 

To assess the efficiency and effectiveness of the algorithms we designed two different sets of experiments.
In the first experiment, we considered instances of \p\ derived from the real-world application data, with randomly generated processing times and varying sizes of the job sets, considering only the shop layout and machines flexibility as in the real-world scenario. 
In the second experiment, we extend the original flow line layout by introducing more flexibility and modifying the relationship between machines and operations.

After describing in more detail the randomly generated instances, we first present the results of some preparatory tests that helped to select and/or discard some variants of the proposed heuristics and then discuss the results of an extensive computational campaign.

\subsection{Instance  description}\label{sec:instances}

In this section, we describe the characteristics of the two sets of instances (hereafter named ``Experiment Set'' 1 and 2) on which the algorithms proposed above have been tested.

\begin{description}
    \item[Experiment Set 1] 
The first set of experiments aims at accurately reproducing the real-world application. So,
for this experiment, we created random test instances following the procedure outlined below. 
We maintained a fixed configuration of five machines with specific inter-stage flexibility constraints, as illustrated in Figure \ref{fig:flowline}. 
To generate these instances, we first set the number of jobs, denoted as $n$, that we wanted to work with. 
Subsequently, we generated the processing times for each of the 16 operations for each job on every possible machine as follows.

In the real-world scenario, operations that can be processed by more than one machine are generally shorter (2 to 14 minutes) than operations that can only be processed by one machine (10 to 28 minutes). 
The latter operations are in fact representative of aggregated operations sets and therefore more time-consuming.
In every instance, a specific operation will represent a task with a certain characteristic, e.g.\ inserting insulation material or fixing a vapor barrier, thus taking short, medium or long time for all jobs. Additionally, each machine has its own features and can be more or less efficient on each operation.
Therefore, for each shiftable operation $i$ and for each machine $M_k$ that can process it, 
we first select an interval \([L_{ik}, H_{ik}]\) by choosing \(L_{ik}\) randomly from a uniform distribution over \([2, 12]\), and then selecting \(H_{ik}\) randomly from a uniform distribution over \([L_{ik}, 14]\).
After settling the characteristics of each pair of operation $i$ and machine $M_k$ this way, we choose the time required to process $o_{ij}$ on machine $M_k$ randomly uniform from $[L_{ik}, H_{ik}]$ for all jobs $j$.
The same approach is used for operations 1, 2, 10, 15, and 16 (i.e., those that can be executed by a single machine $M_k$ only) with $L_{ik}$ chosen from $[10,25]$ and $H_{ik}$ from $[L_{ik}, 28]$.

\item[Experiment Set 2]
In the second set of experiments we extend the original layout of the flow line by adding flexibility to machine $M_4$ and modify the relationship between machines and operations. 
In particular, we assume that for each job the operations $3,4,5,6$ can be executed by machines $M_2$ and $M_3$, operations $8,9,10,11$ by machines $M_3$ and $M_4$, and operations  $13,14,15$ by machines $M_4$ and $M_5$. 
Operations $1,2,7,12$, and $16$ can only be executed by one specific machine, namely machine $M_1,M_2,M_3,M_4$, and $M_5$, respectively. 
Thus, we extend the phases of flexible transition between two machines from $2$ to $3$.
In this setting, the number of feasible assignment modes $|\mathcal{A}_j|$ increases from $40$ to $5\cdot 5\cdot 4=100$. 
This scenario is illustrated in~Figure~\ref{fig:new-layout}. 
The generation of processing times follows the same scheme as for Experiment Set 1. 
\begin{figure}[htbp]
   \centering
       \includegraphics[width=0.8\textwidth]{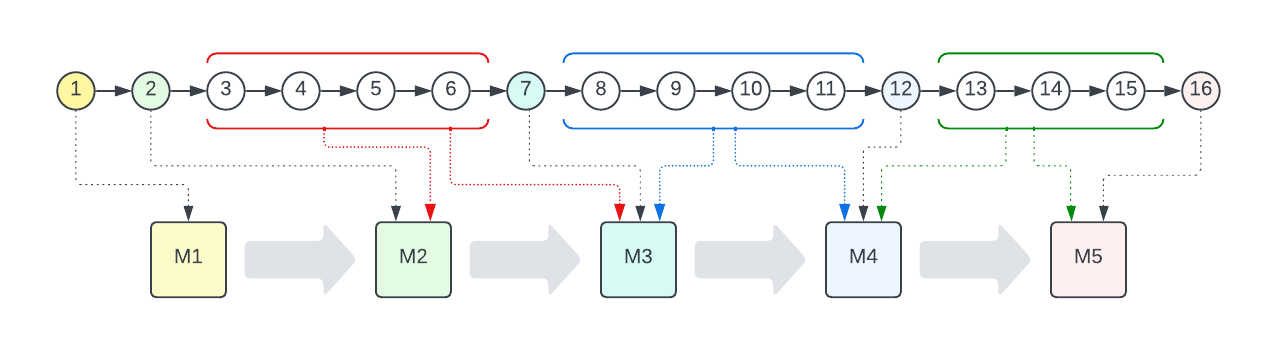}
      \caption{Modified operation-machine relationship of Experiment 2.}\label{fig:new-layout}
\end{figure}
\end{description}

For both sets of experiments, the number $n$ of jobs varies from $5$ to $70$ with stepsize $5$.
For each value of $n$, $10$ instances are generated according to the above described specifics.

In the following sections, for each of the above mentioned 28 classes of instances, we report average values of performance indicators (computation times and objective function values) over the 10 instances of the class.

\subsection{Preparatory experiments}\label{sec:pre-tests}
In order to identify the most promising variants of the heuristics proposed in Sections \ref{sec:heur_ilp} and~\ref{sec:constrheu} and dispose of the less promising algorithmic ideas, a series of experiments have been carried out. 
They provided valuable insights and allowed us to refine our approach by identifying and discarding ineffective strategies.
The best performing variants of these algorithms are then analyzed in a series of more in-depth experiments reported in Section~\ref{sec:comparison_heur}.
Our findings show that: 
\begin{itemize}
\item Firstly, for both the Insertion Heuristic and the RS-OA algorithm (both taking an arbitrary initial job sequence as an input), testing a large number of different initial sequences does not yield significant benefits.
\item Secondly, the four matheuristics presented in Section~\ref{sec:asgnfirst} (i.e., the *-OS heuristics) consistently underperform compared to the Insertion Heuristic and, as it will be clearer below, to all the other methods.
\end{itemize}
In Tables~\ref{tab:single-multi-Insertion} and  \ref{tab:pos_expl_fix_x_random}, we report the performance of the Insertion Heuristic and the matheuristic approach RS-OA, respectively, when running these algorithms once or several times with different starting conditions.
We measure the algorithms' behavior in terms of computation times (column ``Time") and objective function values (column ``Obj.") on instances from Experiment Set 1.  

In particular, Table~\ref{tab:single-multi-Insertion} returns the best result  (minimum objective value) obtained by running the heuristic 50 times, each time giving as an input a new randomly-generated initial sequence, and compare it against the result output from a single shot of the algorithm, i.e., with a single random starting sequence.

\begin{table}[htb]
        \centering
        \begin{tabular}{ccccc}
        \toprule
        &   \multicolumn{2}{c}{Single-start}   &  \multicolumn{2}{c}{Multi-start}  \\  \hline
    $n$	 & 	Time	& 	Obj. 	 & 	{\bf Total} time	&	{\bf Best} Obj.	\\ 
    \midrule
	5	& 0.01	&   408.7 	&	0.29	&	404.5  \\
	10	& 0.03	&	 674.7	&	1.73	&	666.2	\\
	15	& 0.08	&	 995.7	&	4.87	&	984.1	\\
	20	& 0.18	&	1237.6	&	11.76	&	1224.0	\\
	25	& 0.36	&	1596.4	&	22.41	&	1578.4	\\
    \bottomrule
        \end{tabular}
        \caption{Single vs. multi-start (50 input sequences) results for the Insertion Heuristic.}
        \label{tab:single-multi-Insertion}
    \end{table}

Similarly, Table \ref{tab:pos_expl_fix_x_random} compares the results 
for the matheuristic approach RS-OA. Recall this approach fixes the job sequencing and then uses MIP1 to optimally solve the assignment of operations to machines. 
Also in this case, we compare the objective value obtained starting from a single randomly-generated sequence against the best result out of 50 distinct solutions, obtained starting from different sequences. 

\begin{table}[htb]
        \centering
        \begin{tabular}{ccccc}
        \toprule
        &   \multicolumn{2}{c}{Single-start}   &  \multicolumn{2}{c}{Multi-start}  \\  \hline
$n$	 & 	Time	& 	Obj. 	 & 	{\bf Total} time	&	{\bf Best} Obj.	\\ 
\midrule
	5	&	0.04	&	384.3	&	2.09	&	370.5	\\
	10	&	0.16	&	652.1	&	9.99	&	639.0	\\
	15	&	0.87	&	931.6	&	52.42	&	914.8	\\
	20	&	4.44	&	1162.4	&	247.94	&	1149.9	\\
	25	&	17.03	&	1489.9	&	725.14	&	1474.0	\\
    \bottomrule
        \end{tabular}
        \caption{Single vs. multi-start (50 input sequences) results for the Random  Sequencing with Optimal  Assignment (RS-OA) matheuristic.}
    \label{tab:pos_expl_fix_x_random}
    \end{table}

These results clearly show how the initial job-sequence decision 
has a negligible impact on the value of the objective function returned by both the Insertion Heuristic and the RS-OA algorithms. 
The improvement when these heuristics are run on fifty initial sequences (and the best solution is returned) is marginal, with only an enhancement of at most 
1\% for the Insertion Heuristic and 4\% for the RS-OA in the objective function compared to a significantly larger (between $40$ and $60$ times for both methods) computation time. 
In addition, the ratio between the multi-start and the single-start objective values remains essentially constant as $n$ increases, so we deemed not necessary to report the results for larger values of $n$.
Based on the above results, 
we limit the remaining analysis of the performance for the RS-OA and the Insertion Heuristic to their single-start variants, i.e., when only one starting random job-sequence is used.

\medskip
Hereafter, we show that the matheuristic algorithms which first determine the operation assignment and then compute an optimal job sequence (see Section~\ref{sec:asgnfirst}), consistently perform worse than alternative methods. 
\begin{table}[htb]    
        \footnotesize
        \centering
        \begin{tabular}{crrrrrrrrrrr}
        \toprule
        &   \multicolumn{2}{c}{RA-OS}   &  \multicolumn{2}{c}{LA-OS} & \multicolumn{2}{c}{IA-OS} & \multicolumn{3}{c}{MA-OS} & \multicolumn{2}{c}{Insertion Heur.} 
        \\  \hline
    $n$	 &  Time 	& 	Obj.	  & 	Time	&	Obj.  & 	Time	&	Obj. & 	 Time	&	Obj.  & Opt. & Time & Obj. 
    \\ \midrule
    5  &         0.04 &    501.0 &  0.06 &   473.0 &       0.04 &  501.3 &   0.04 &  425.0  &  100 \% & \textbf{0.01} & \textbf{ 408.7} 
    \\
    10 &         0.12 &    888.6 &  0.17 &   876.4 &       0.15 &  847.1 &   0.19 &  769.4  &  100 \% & \textbf{0.03} & \textbf{ 674.7} 
    \\
    15 &         0.29 &   1286.4 &  0.31 &  1259.5 &       1.36 & 1210.0 &   0.34 & 1190.3  &  100 \% & \textbf{0.08} & \textbf{ 995.7} 
    \\
    20 &         0.47 &   1658.8 &  0.54 &  1676.5 &       0.90 & 1587.7 &  40.48 & 1385.0  &   90 \% & \textbf{0.18} & \textbf{1237.6} 
    \\
    25 &         0.76 &   2066.3 &  0.75 &  2195.5 &       1.55 & 2141.0 &  43.33 & 2021.3  &   90 \% & \textbf{0.36} & \textbf{1596.4} 
    \\
    30 &         1.31 &   2575.2 & 22.44 &  2168.1 &      20.13 & 2168.2 &   1.05 & 2261.9  &  100 \% & \textbf{0.60} & \textbf{1885.9} 
    \\
    35 &         1.37 &   2766.3 &  5.12 &  2669.7 &       3.71 & 2670.6 &   1.50 & 2959.7  &  100 \% & \textbf{0.93} & \textbf{2187.8} 
    \\
    40 &         2.06 &   3172.8 &  8.90 &  2876.1 &       7.76 & 2876.1 & 125.52 & 2717.1  &   60 \% & \textbf{1.34} & \textbf{2404.5} 
    \\
    45 &         2.29 &   3799.5 & 11.05 &  3525.2 &       9.43 & 3525.3 &  63.62 & 3447.9  &   90 \% & \textbf{1.97} & \textbf{2818.5} 
    \\
    50 &         3.80 &   4116.8 & 12.95 &  3648.0 &      11.83 & 3648.3 &  66.07 & 3343.1  &   80 \% & \textbf{2.65} & \textbf{2942.5} 
    \\
    55 &         6.45 &   4633.7 &  8.07 &  4115.7 &       9.28 & 4397.4 & 125.17 & 3914.7  &   60 \% & \textbf{3.74} & \textbf{3373.3} 
    \\
    60 &         6.72 &   4981.2 & 12.46 &  4753.3 &      19.48 & 4681.5 &  99.56 & 4255.5  &   70 \% & \textbf{4.43} & \textbf{3604.5} 
    \\
    65 & {5.95}&   5260.9 & 23.61 &  4756.6 &      17.60 & 4756.6 & 118.02 & 4790.0  &   70 \% & \textbf{5.95} & \textbf{3898.9} 
    \\
    70 & {\bf 6.04}&   5740.3 & 38.32 &  4950.0 &      34.51 & 4949.4 &  71.67 & 5210.8  &   80 \% &         7.45  & \textbf{4094.0} 
    \\
    \bottomrule
        \end{tabular}
        \caption{Comparison of matheuristics based on fixing assignment first.}\label{tab:fix_assignment}
    \end{table}

Table~\ref{tab:fix_assignment} presents the results for the first type of heuristics 
(specifically the RA-OS, IA-OS, LA-OS, and MA-OS heuristics) for the instances in Experiment Set 1,  together with the results for the Insertion Heuristic as a comparison.
The computation time comprises both the time spent in determining the fixed operation assignment and the time spent on computing the optimal job sequence by solving the restricted instance of MIP1.
However the contribution of the first phase is almost negligible for the RA-OS, LA-OS and MA-OS matheuristics. 
Interestingly enough, once the operation assignment is done, the remaining sequencing problem can be optimally solved within a short amount of time for almost all instances by all the algorithms but MA-OS. 
In column ``Opt.", we report the fraction of instances in which MA-OS optimally solves the sequencing step  within the time limit of 1800 seconds, for an operation assignment given in the first phase. 
In general, the CPU time required by this final sequencing step is subject to large variations across instances of the same size. 
Even more interestingly, larger CPU times seem to be associated with better assignment choices, as can be observed by looking at the results of the MA-OS matheuristic. 

The performance of these four algorithms can be compared to that of the Insertion Heuristic which is reported in the last columns of Table~\ref{tab:fix_assignment} (and also in Table~\ref{tab:Exp1_results}).
Insertion Heuristic clearly dominates RA-OS, LA-OS, IA-OS, and MA-OS in terms of both efficiency, measured by computation time in seconds, and effectiveness, evaluated by the average makespan, i.e., the solution objective value. 
The only exception being the computation time of RA-OS on the instances with $n=70$ which is slightly less than that of the Insertion Heuristic but together with a dramatic worsening in the objective function values. 

Based on the above results, we excluded RA-OS, IA-OS, LA-OS, and MA-OS heuristics from further analyses.

\subsection{Comparison of the proposed solution approaches}\label{sec:comparison_heur}

We analyze the performance of the heuristic algorithms presented in Sections~\ref{sec:seqfirst},  \ref{sec:holistic}, and~\ref{sec:constrheu}, against the MIP models illustrated in Section~\ref{sec:mip} under a time limit of 1800 seconds. 

\begin{table}[htbp]
    \centering
        \begin{tabular}{rrrrr}
\toprule
      & \multicolumn{3}{c}{Objective value after} & Optimally \\
{$n$} &  60 s &  300 s &  1800 s & solved \\
\midrule
5  &            368.8 &             368.8 &          368.8 & 100\%  \\
10 &            625.2 &             624.5 &          624.5 & 100\%  \\
15 &            905.0 &             899.6 &          897.4 &  20\%  \\
20 &           1140.9 &            1132.9 &         1127.6 &   0\%  \\
25 &           1477.3 &            1459.2 &         1451.0 &   0\%  \\
30 &           1779.3 &            1719.1 &         1703.0 &   0\%  \\
35 &           2051.0 &            2002.9 &         1977.7 &   0\%  \\
40 &           2266.9 &            2221.1 &         2189.1 &   0\%  \\
45 &           2639.6 &            2585.2 &         2548.3 &   0\%  \\
50 &           2878.6 &            2818.2 &         2768.4 &   0\%  \\
55 &           3147.3 &            3107.9 &         3032.8 &   0\%  \\
60 &           3454.5 &            3422.0 &         3372.4 &   0\%  \\
65 &           3766.8 &            3700.9 &         3587.5 &   0\%  \\
70 &           4084.0 &            4009.4 &         3862.7 &   0\%  \\
\bottomrule
\end{tabular}
    \caption{{\bf Experiments set 1}. Results for the positional implicit model MIP2 after 60, 300 and 1800 seconds.}
    \label{tab:MIP-Exp1}
\end{table}

Table~\ref{tab:MIP-Exp1} reports the average objective values obtained by MIP2, the best performing MIP model according to \cite{us_ods2023}, after one, five and thirty minutes of computation, for the 14 (pseudo-randomly generated) classes of instances with 5 up to 70 jobs described in Section~\ref{sec:instances}. 
For instance classes with at least 15 jobs, the solver is not able to certify optimality of found solutions within the allotted time. 
As it is expected, larger computation times allow to obtain better incumbent solutions. 
However, larger instances prove to be quite hard for the solver, as half an hour of computation permits an improvement of solutions only by approximately 5\%. 

\begin{sidewaystable}
\centering			
{\footnotesize		
\begin{tabular}{c rr rr rr rr rr rr rr}	
\toprule
	&	\multicolumn{8}{c }{SFT}	&		&		&		&		&		&		 \\
	&	\multicolumn{2}{c}{SFT(1,0.66)}					&	\multicolumn{2}{c}{SFT(1,0.51)}			&	\multicolumn{2}{c}{SFT($n$,0.66)}	&	\multicolumn{2}{c}{SFT($n$,0.51)}			&	\multicolumn{2}{c }{RS-OA}			&	\multicolumn{2}{c }{MS-OA}			&	\multicolumn{2}{c }{Insertion Heur.}			\\ \hline
$n$	&	Time	&	Obj	&	Time	&	Obj	&	Time	&	Obj	&	Time	&	Obj	&	Time	&	Obj	&	Time	&	Obj	&	Time	&	Obj	\\
5	&	0.23	&	{\bf 374.0}	&	0.15	&	391.3	&	0.19	&	375.1	&	{\bf 0.14}	&	392.0	&	{\bf 0.04}	&	384.3	&	0.08	&	375.3	&	0.01	&	408.7	\\
10	&	6.23	&	{\bf 632.4}	&	0.70	&	667.6	&	7.28	&	635.3	&	0.60	&	669.0	&	{\bf 0.16}	&	652.1	&	0.40	&	642.0	   &	0.03	&	674.7	\\
15	&	171.6	&	{\bf 906.3}	&	9.26	&	967.0	&	137.84	&	910.4	&	1.54	&	973.2	&	{\bf 0.87}	&	931.6	&	1.62	&	920.7	   &	0.08	&	995.7	\\
20	&	206.93	&	{\bf 1137.4}	&	5.34	&	1202.1	&	93.10	&	1147.2	&	2.86	&	1226.8	&	{\bf 4.44}	&	1162.4	&	28.93	&	1152.7	   &	0.18	&	1237.0	\\
25	&	276.40	&	{\bf 1458.4}	&	20.65	&	1549.1	&	285.74	&	1459.8	&	12.09	&	1579.8	&	{\bf 17.03}	&	1489.9	&	40.49	&	1481.3	   &	0.36	&	1596.0	\\
30	&	287.83	&	{\bf 1708.8}	&	10.35	&	1829.9	&	289.75	&	1710.1	&	7.05	&	1870.9	&	50.01	&	1751.7	&	{\bf 47.03}	&	1736.2	   &	0.60	&	1885.90	\\
35	&	259.70	&	1988.8	&	17.91	&	2105.6	&	257.98	&	{\bf 1986.8}	&	9.36	&	2171.0	&	{\bf 67.62}	&	2025.8	&	68.64	&	2016.5	   &	0.93	&	2187.8	\\
40	&	314.08	&	2193.6	&	26.53	&	2303.3	&	312.81	&	{\bf 2192.3}	&	12.75	&	2340.4	&	{\bf 84.98}	&	2225.5	&	157.58	&	2216.4	   &	1.34	&	2404.5	\\
45	&	291.99	&	{\bf 2569.2}	&	45.96	&	2720.7	&	272.02	&	2579.9	&	17.20	&	2801.0	&	{\bf 75.81}	&	2583.4	&	144.68	&	2572.1	   &	1.97	&	2818.5	\\
50	&	322.91	&	2777.2	&	42.05	&	2904.4	&	264.01	&	2788.3	&	21.28	&	2963.9	&	{\bf 90.64}	&	2788.5	&	135.86	&	{\bf 2776.5}	   &	2.65	&	2942.5	\\
55	&	305.33	&	3040.5	&	51.88	&	3159.9	&	302.70	&	3029.9	&	28.21	&	3204.2	&	{\bf 106.69}	&	3048.6	&	186.98	&	{\bf 3023.0}	&	3.74	&	3373.3	\\
60	&	320.97	&	3407.6	&	80.49	&	3525.6	&	310.02	&	3400.1	&	33.87	&	3604.6	&	{\bf 115.89}	&	3364.8	&	178.53  &	{\bf 3354.6}	&	4.43	&	3604.5	\\
65	&	325.41	&	3625.1	&	117.76	&	3795.5	&	337.37	&	3623.7	&	39.10	&	3875.1	&	{\bf 117.04}	&	3596.2	&	204.20	&	{\bf 3579.9}	&	5.95	&	3898.9	\\
70	&	335.43	&	3864.4	&	132.24	&	4060.2	&	311.65	&	3872.2	&	45.84	&	4137.3	&	{\bf 105.46}	&	3858.2	&	178.23	&	{\bf 3850.5}	&	7.45	&	4094.0	\\
\bottomrule																				
\end{tabular}				
}
\caption{{\bf Experiment Set 1}. Computational results.}
\label{tab:Exp1_results}.
\end{sidewaystable}

Table~\ref{tab:Exp1_results} reports the same performance indicators, namely, computation times (col.\ ``Time") and average objective values (col.\ ``Obj.") obtained by the proposed heuristic algorithms on the same 14 classes of instances.
 As before, values are average over 10 instances for each class. 
 Regarding the SFT algorithm, we focus on examining two threshold values: \(\phi = 0.66\) (high threshold) and \(\phi = 0.51\) (low threshold), as well as the extreme cases where only one variable or all variables are fixable, i.e., \(r = 1\) and \(r = n\). 
We consider the resulting four configurations to be both important and representative of the algorithm's behavior.

Comparing the objective values obtained by the best performing heuristic against the corresponding results of the MIP2 model processed with the Gurobi solver---after 1800 seconds of computation---we may note a slightly superior behavior of the latter method for instances with up to 55 jobs. 
In this case, the little improvement over the best heuristics ranges from 0.1\% to 1.4\% (with the highest improvement seen at $n=5$). However, the limited prevalence of the MIP model comes at the expense of significantly increased computation time.
For example, for $n=50$, SFT$(1, 0.66)$ only takes 322.91 seconds to reach an only slightly worse solution than the MIP approach with the 1800 seconds time limit.
In general, for small instances, if computation time is not an issue, the MIP model still provides the best accuracy. On the other hand, when time has an impact, SFT(1, 0.66) is a faster option that offers results close to the best, especially for smaller instances.
When dealing with large instances, the MIP model is outperformed by the best heuristics even at the 1800 seconds mark, still consistently delivering solutions whose quality is close to (between 0.2\% and 0.5\% worse than) the best. In addition, half an hour computation for the MIP model must be confronted to the time taken by other methods, like MS-OA and RS-OA, achieving similar or slightly better results in much less time (see for instance, MS-OA solving $n=70$ in 178.23 seconds, with a slightly better objective value than the MIP's).
\emph{A fortiori}, restricting our attention to results obtained under a computation time limited to around 300 seconds, we may note that the MIP at the 300-second mark is dominated by other methods, like SFT(1,0.66) and MS-OA,  making them the preferred choice for both accuracy and time efficiency, especially on larger problems.

Summarizing, for small instances up to 50 jobs, the solver applied to MIP2 obtains slightly better results than the heuristics, in terms of solution values, but it is significantly slower. Sacrificing just a bit in terms of solution quality, SFT(1, 0.66) provides good solutions in a much faster way. For large instances, the MIP model still performs well in terms of objective values, but it is not competitive anymore against other methods, as MS-OA offers a better balance between time and solution quality. In conclusion, other methods, especially MS-OA and SFT(1, 0.66), strike a better tradeoff between accuracy and time.

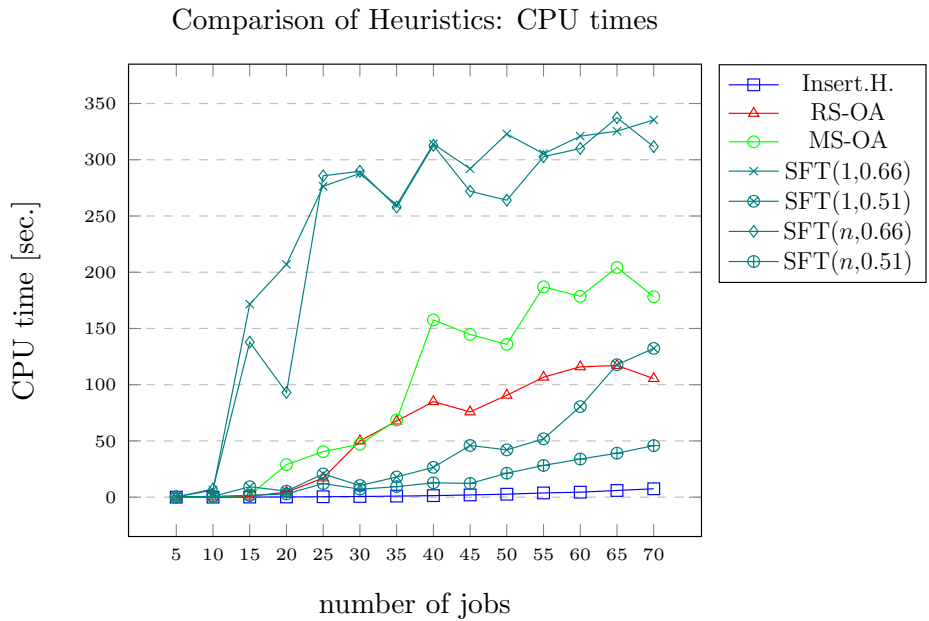
\begin{figure}[htbp]
    \centering
\begin{tikzpicture}[scale=1.1] 
\begin{axis}
[
    title={\footnotesize Comparison of Heuristics: CPU times},
    xlabel={\footnotesize number of jobs},
    ylabel={\footnotesize CPU time [sec.]},
    xmin=5, xmax=70,
    ymin=0, ymax=350,
    enlargelimits=0.1,
    symbolic x coords={5,10,15,20,25,30,35,40,45,50,55,60,65,70},
    xtick={5,10,15,20,25,30,35,40,45,50,55,60,65,70},
    ytick={0,50,100,150,200,250,300,350},
    legend pos=outer north east,
    ymajorgrids=true,
    grid style=dashed,
    legend style={nodes={scale=0.7}}
]
\addplot[
    color=blue, 
    mark=square,
    ]
    coordinates {
    (5,0.01)(10,0.03)(15,0.08)(20, 0.18)(25,0.36)(30, 0.60)(35, 0.93)(40,1.34)(45,1.97)(50,2.65)(55,3.74)(60,4.43)(65,5.95)(70,7.45)
    };
\addplot[
    color=red, 
    mark=triangle,
    ]
    coordinates {
    (5,0.04)(10,0.16)(15,0.87)(20,4.44)(25,17.03)(30,50.01)(35,67.62)(40,84.98)(45,75.81)(50,90.64)(55,106.69)(60,115.89)(65,117.04)(70,105.46)
    };
\addplot[
    color=green, 
    mark=o,
    ]
    coordinates {
    (5,0.08)(10,0.40)(15,1.62)(20,28.93)(25,40.49)(30,47.03)(35,68.64)(40,157.58)(45,144.68)(50,135.86)(55,186.98)(60,178.53)(65,204.20)(70,178.23)
    };
\addplot[
    color=teal, 
    mark=x,
    ]
    coordinates {
    (5,0.23)(10,6.23)(15,171.60)(20,206.93)(25,276.40)(30,287.83)(35,259.70)(40,314.08)(45,291.99)(50,322.91)(55,305.33)(60,320.97)(65,325.41)(70,335.43)
    };
\addplot[
    color=teal, 
    mark=otimes,
    ]
    coordinates {
    (5,0.15)(10,0.70)(15,9.26)(20,5.34)(25,20.65)(30,10.35)(35,17.91)(40,26.53)(45,45.96)(50,42.05)(55,51.88)(60,80.49)(65,117.76)(70,132.24)
    };
\addplot[
    color=teal, 
    mark=diamond,
    ]
    coordinates {
    (5,0.19)(10,7.28)(15,137.84)(20,93.10)(25,285.74)(30,289.75)(35,257.98)(40,312.81)(45,272.02)(50,264.01)(55,302.70)(60,310.02)(65,337.37)(70,311.65)
    };
\addplot[
    color=teal, 
    mark=oplus,
    ]
    coordinates {
    (5,0.14)(10,0.60)(15,1.54)(20,2.86)(25,12.09)(30,7.05)(35,9.36)(40,12.75)(45,12.20)(50,21.28)(55,28.21)(60,33.87)(65,39.10)(70,45.84)
    };
\addlegendentry{Insert.H.}
\addlegendentry{RS-OA}
\addlegendentry{MS-OA}
\addlegendentry{SFT(1,0.66)}
\addlegendentry{SFT(1,0.51)}
\addlegendentry{SFT($n$,0.66)}
\addlegendentry{SFT($n$,0.51)}
\end{axis}
\end{tikzpicture}
\par\vskip12pt
\caption{{\bf Experiment Set 1}. CPU times comparison of Heuristics (plots of data in Table~\ref{tab:Exp1_results}).}
    \label{fig:Exp1-CPUplot}
\end{figure}
\begin{figure}\centering
\begin{tikzpicture}[scale=1.1] 
\begin{axis}
[
    title={\footnotesize Comparison of Heuristics: difference from best objective},
    xlabel={\footnotesize number of jobs},
    ylabel={\footnotesize difference in average objectives},
    xmin=5, xmax=70,
    ymin=0, ymax=350,
    enlargelimits=0.1,
    symbolic x coords={5,10,15,20,25,30,35,40,45,50,55,60,65,70},
    xtick={5,10,15,20,25,30,35,40,45,50,55,60,65,70},
    ytick={0,50,100,150, 200, 250, 300,350},
    legend pos=outer north east,
    ymajorgrids=true,
    grid style=dashed,
    legend style={nodes={scale=0.7}}
]
\addplot[
    color=blue, 
    mark=square,
    ]
    coordinates {
    (5,34.7)(10,42.3)(15,89.4)(20,100.2)(25,138.0)(30,177.1)(35,201.0)(40,212.2)(45,249.3)(50,166)
    (55,350.3)(60,249.9)(65,319)(70,243.5)
    };
\addplot[
    color=red, 
    mark=triangle,
    ]
    coordinates {
    (5,10.3)(10,19.7)(15,25.3)(20,25.0)(25,31.5)(30,42.9)(35,39.0)(40,33.2)(45,14.2)    (50,12.0)
    (55,25.6)(60,10.2)(65,16.3)(70,7.7)
    };
\addplot[
    color=green, 
    mark=o,
    ]
    coordinates {
    (5,1.3)(10,9.6)(15,14.4)(20,15.3)(25,22.9)(30,27.4)(35,29.7)(40,24.1)(45,2.9)(50,0)
    (55,0)(60,0)(65,0)(70,0)
    };
\addplot[
    color=teal, 
    mark=x,
    ]
    coordinates {
    (5,0)(10,0)(15,0)(20,0)(25,0)(30,0)(35,2)(40,1.3)(45,0)    (50,0.7)
    (55,17.5)(60,53.0)(65,45.2)(70,13.9)
    };
\addplot[
    color=teal, 
    mark=otimes,
    ]
    coordinates {
    (5,17.3)(10,35.2)(15,60.7)(20,64.7)(25,90.7)(30,121.1)(35,118.8)(40,111.0)(45,151.5)    (50,127.9)
    (55,136.9)(60,171.0)(65,215.6)(70,209.7)
    };
\addplot[
    color=teal, 
    mark=diamond,
    ]
    coordinates {
    (5,1.1)(10,2.9)(15,4.1)(20,9.8)(25,1.4)(30,1.3) (35,0)(40,0)(45,10.7)    (50,11.8)
    (55,6.9)(60,45.5)(65,43.8)(70,21.7)
    };
\addplot[
    color=teal, 
    mark=oplus,
    ]
    coordinates {
    (5,18)(10,36.6)(15,66.9)(20,89.4)(25,121.4)(30,162.1)(35,184.2)(40,148.1)(45,231.8)    (50,187.4)
    (55,181.2)(60,250)(65,295.2)(70,286.8)
    };
\addlegendentry{Insert.H.}
\addlegendentry{RS-OA}
\addlegendentry{MS-OA}
\addlegendentry{SFT(1,0.66)} 
\addlegendentry{SFT(1,0.51)} 
\addlegendentry{SFT($n$,0.66)} 
\addlegendentry{SFT($n$,0.51)} 
\end{axis}
\end{tikzpicture}
\caption{{\bf Experiment Set 1}. Objective values comparison: Absolute differences in terms of average values of  objective functions from best average objective (plots derived from data in Table~\ref{tab:Exp1_results}).}
\label{fig:Exp1-ObjDiffPlot}
\end{figure}
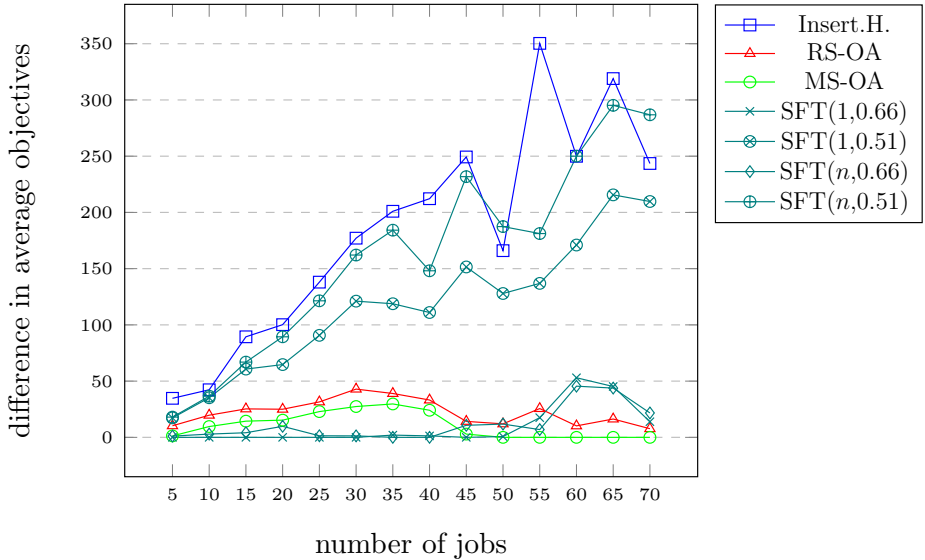

\smallskip
In order to provide more detailed insights into the tradeoff between computation time and objective values for each method, let's dive deeper into the specific performance of the heuristic algorithms (including the Insertion Heuristic), comparing their results on small and large instances. Figures~\ref{fig:Exp1-CPUplot} and~\ref{fig:Exp1-ObjDiffPlot} give a pictorial representation of the data provided in Table~\ref{tab:Exp1_results} which we comment hereafter.

\begin{itemize}
    \item As for the SFT$(r,\phi)$ algorithms, we immediately note how a low threshold makes the procedure faster but produces worse results, in terms of solution quality. 
    On average, SFT(1,0.66) is 5.41\% more accurate than SFT(1,0.51), and SFT($n$,0.66) is 6.71\% more accurate than SFT($n$,0.51) in terms of objective values across the tested instance sizes. 
    Greater solution quality with lower threshold comes at the cost of increased computation times: SFT(1,0.51) and SFT($n$,0.51) are approximately 13 and 26 faster than SFT(1,0.66) and SFT($n$,0.66), respectively, on average across the tested instance sizes. This significant decrease in computation time reflects the trade-off between solution accuracy and processing speed.
    \item Insertion Heuristic is---as it can be expected---the fastest method by far, with running times much lower than the others. 
    However, the objective values are consistently worse than those obtained with other methods. The reduction in terms of computation times is substantial, especially for larger instances. 
    Insertion Heuristic is from 75.0\% to 97.02\%  faster than the fastest one among the other methods, with an average reduction in terms of computation time of 88.49\%. 
    Unfortunately, this method provides solution with larger objective values. Compared to the best results obtained by the companion algorithms, we measure an accuracy loss which is between 5.97\% and  10.37\% with an average of  8.79\% increase in the objective values. 
    \item An overall qualitative comparison is illustrated by Figures~\ref{fig:Exp1-CPUplot} and~\ref{fig:Exp1-ObjDiffPlot} in which a comparison among the solutions' quality of the various approaches is made more evident reporting the absolute differences from the best result across every instance class. 
    We may distinguish three main groups of algorithms with increasing computation times corresponding to larger quality (i.e., smaller average objective values) of the provided solutions. SFT(1,0.66) and SFT($n$,0.66) are the most time-expensive, while RS-OA and MS-OA are approximately twice faster. 
    The SFT algorithms with a low threshold, i.e., SFT(1,0.51) and SFT($n$,0.51), are consistently the fastest algorithms across all the instances, obviously, disregarding the Insertion Heuristic.
    Not surprisingly, the threshold value $\phi$ has a noticeable impact on the speed, as it allows to fix a larger resp.\ smaller number of variables in MIP2 at the cost of a decreased resp.\ increased number of iterations and hence in the accuracy of the remaining models. 
    The situation is somehow overturned if we consider quality of solutions (see Figure~\ref{fig:Exp1-ObjDiffPlot} 
  in which, for each algorithm, we plot the difference between the average objective values of the algorithm and of the best performing model). Again SFT(1,0.66) and SFT($n$,0.66) show similar results and are the most accurate algorithms for small instances up to 40 jobs while kindred RS-OA and MS-OA become the winners for larger instances. 
    Instead, the SFT methods with $\phi=0.51$ and the Insertion Heuristic are the worst performing algorithm. 
\end{itemize}
In conclusion, RS-OA and MS-OA strike a balance between time and objective value, performing relatively well on both measures, with MS-OA showing competitive objective values for larger instances ($n\geq 45$). The time required for all methods obviously increases with $n$, but some methods---like RS-OA or SFT algorithms with a low threshold ($\phi=0.51$)---scale better than others.


\begin{table}[htbp]
    \centering
 \begin{tabular}{lrrrr}
    \toprule
 & \multicolumn{3}{c}{Objective value after} & Optimally \\
{$n$} &  60 s &  300 s &  1800 s & solved \\
\midrule
    5  &            349.3 &             349.3 &              349.3 &   100\% \\
    10 &            566.5 &             564.7 &              564.1 &    30\% \\
    15 &            797.3 &             792.9 &              790.3 &     0\% \\
    20 &           1018.0 &            1008.0 &             1003.4 &     0\% \\
    25 &           1326.0 &            1300.9 &             1290.5 &     0\% \\
    30 &           1444.5 &            1414.9 &             1402.6 &     0\% \\
    35 &           1753.2 &            1731.5 &             1697.6 &     0\% \\
    40 &           1979.9 &            1961.1 &             1905.0 &     0\% \\
    45 &           2260.5 &            2227.8 &             2171.0 &     0\% \\
    50 &           2556.4 &            2473.5 &             2407.8 &     0\% \\
    55 &           2795.6 &            2745.1 &             2683.4 &     0\% \\
    60 &           3320.9 &            3072.4 &             2993.4 &     0\% \\
    65 &           3534.3 &            3239.6 &             3162.2 &     0\% \\
    70 &           3476.3 &            3394.1 &             3316.3 &     0\% \\
    \bottomrule
\end{tabular}
\caption{{\bf Experiments set 2}. Results for the positional implicit model MIP2 after 60, 300 and 1800 seconds.}\label{table:pos_impl_2}
\end{table}

\begin{sidewaystable}
\centering			
{\footnotesize		
\begin{tabular}{c rr rr rr rr rr rr rr}	
\toprule
	&	\multicolumn{8}{c}{SFT}	&		&		&		&		&		&		 \\
	&	\multicolumn{2}{c}{SFT$(1,0.66)$}					&	\multicolumn{2}{c }{SFT$(1,0.51)$}			&	\multicolumn{2}{c }{SFT$(n,0.66)$}	&	\multicolumn{2}{c }{SFT$(n,0.51)$}			&	\multicolumn{2}{c }{RS-OA}			&	\multicolumn{2}{c }{MS-OA}			&	\multicolumn{2}{c }{Insertion Heur.}			\\ \hline
$n$	&	Time	&	Obj	&	Time	&	Obj	&	Time	&	Obj	&	Time	&	Obj	&	Time	&	Obj	&	Time	&	Obj	&	Time	&	Obj	\\
5  &       0.82 &    {\bf 351.8} &       0.49 &   357.3 &       0.72 &  {\bf   351.8} &    0.48 &   361.1 &   0.04 &          362.2 &       0.14 &          354.8 &   {\bf 0.01} &          368.2 \\ 
10 &     134.37 &    {\bf 564.4} &      51.14 &   577.5 &     115.66 &          566.8 &   31.98 &   583.6 &   0.23 &          584.8 &       0.95 &          577.3 &   {\bf 0.03} &          615.5 \\ 
15 &     295.27 &    {\bf 793.0} &     111.11 &   813.2 &     275.88 &          793.3 &    9.14 &   826.5 &   1.47 &          813.6 &       6.74 &          808.1 &   {\bf 0.09} &          864.8 \\ 
20 &     307.49 &         1006.0 &      82.49 &  1036.4 &     295.93 &  {\bf  1005.8} &   40.48 &  1046.2 &   4.37 &         1036.4 &      50.81 &         1026.1 &   {\bf 0.18} &         1107.8 \\ 
25 &     313.27 &         1300.9 &     195.33 &  1327.4 &     311.67 &  {\bf  1296.8} &   52.48 &  1354.9 &  44.57 &         1319.1 &      71.87 &         1312.5 &   {\bf 0.35} &         1384.1 \\ 
30 &     317.36 &         1422.6 &     142.38 &  1438.1 &     316.67 &  {\bf  1420.4} &  115.78 &  1464.0 & 105.13 &         1437.3 &     166.39 &         1423.5 &   {\bf 0.58} &         1480.9 \\ 
35 &     325.66 &         1724.1 &     128.13 &  1743.6 &     323.50 &         1718.2 &   35.66 &  1790.2 &  96.13 &         1728.9 &     159.02 &  {\bf  1717.2} &   {\bf 0.92} &         1817.7 \\ 
40 &     333.33 &         1934.4 &     204.58 &  1928.7 &     329.93 &         1926.0 &   77.08 &  1982.5 & 115.72 &         1929.8 &     179.92 &  {\bf  1915.8} &   {\bf 1.40} &         2022.3 \\ 
45 &     342.40 &         2233.8 &     257.94 &  2213.8 &     339.14 &         2222.0 &  101.71 &  2260.1 & 115.37 &         2186.4 &     199.71 &  {\bf  2173.7} &   {\bf 2.04} &         2270.0 \\ 
50 &     354.02 &         2468.7 &     226.20 &  2462.9 &     350.92 &         2458.1 &  127.01 &  2498.4 & 120.00 &         2400.9 &     205.78 &  {\bf  2394.5} &   {\bf 2.77} &         2565.0 \\ 
55 &     362.64 &         2754.3 &     297.52 &  2707.1 &     359.82 &         2730.9 &  111.59 &  2770.0 & 117.13 &         2680.4 &     195.55 &  {\bf  2666.4} &   {\bf 3.58} &         2858.5 \\ 
60 &     378.95 &         3042.2 &     369.20 &  2994.9 &     383.46 &         3022.3 &  141.70 &  3094.0 & 111.98 &         2968.9 &     169.34 &  {\bf  2960.8} &   {\bf 4.64} &         3147.4 \\ 
65 &     394.10 &         3199.8 &     422.73 &  3158.0 &     374.72 &         3261.0 &  236.54 &  3205.1 & 119.83 &         3128.8 &     199.77 &  {\bf  3119.2} &   {\bf 5.95} &         3330.1 \\ 
70 &     409.23 &         3366.9 &     553.14 &  3267.5 &     410.52 &         3390.3 &  186.50 &  3345.2 & 123.67 &         3279.7 &     257.40 &  {\bf  3254.1} &   {\bf 8.46} &         3405.4 \\ 
\bottomrule																						
\end{tabular}				
}
\caption{{\bf Experiment Set 2}. Computational results.}
\label{tab:Exp2_res}
\end{sidewaystable}												
\bigskip
The results from the second set of instances (Experiment Set 2) largely confirm the trends observed in Set 1 and are presented below. 

Regarding the results obtained by solving the MIP2 model with Gurobi, similar observations to those made for Experiment Set 1 also apply to the second set.
If we compare the solution quality (objective values) of the MIP against the best performing heuristic, disregarding the computation time, we may note a tiny (between 0 and 1.1\%) advantage of the MIP for instances up to 50 job. 
On the other hand 30 minutes of computation time is required to gain such a minor improvement whereas the best heuristics output their solutions in less than 200 seconds.
For larger instances with $n\geq 55$ jobs, the MIP model is dominated on both indicators by the best heuristic algorithms.

In order to highlight the key differences in the trends of the heuristic algorithms' performance between the two experiments, again we compare both the computation time and objective value across the algorithms for different problem sizes. 
Once more, we may distinguish the same three classes of algorithms characterized by increasing accuracy and increasing computation times. 

Looking at computation times of the algorithms, with respect to the first set of experiments, we note the following
\begin{itemize}
    \item SFT(1,0.51) is still generally faster than SFT(1,0.66), but the differences are less pronounced for Experiment Set 2. 
    Additionally, times for both versions of SFT are significantly higher for smaller instances when compared to Experiment 1 (e.g., 553.14 for $n=70$ in Experiment 2 compared to 132.24 in Experiment 1). 
    With regard to SFT($n$,0.66) and SFT($n$,0.51), we also observe a more consistent increase in CPU times than in the first experiment. 
    SFT($n$,0.51) is again faster than SFT($n$,0.66), but the difference is not as large as for the versions with $r=1$. 
    Even the faster SFT($n$,0.51) exhibits a significant increase in time (e.g., from 48.87 at $n=5$ to 186.5 at $n=70$) in Experiment Set 2. In general, the gap between the two experiments becomes more evident at larger problem sizes.
    \item RS-OA and MS-OA computation times increase at a steadier rate compared to Experiment 1, but the differences in computation time between them and SFT algorithms become less pronounced. 
    Both algorithms remain competitive in the larger problem sizes as well, especially MS-OA.
    \item The Insertion Heuristic obviously remains the fastest in terms of computation time, maintaining the trend from Experiment 1.
\end{itemize}
As for the solutions' objective function values we can observe:
\begin{itemize}
    \item The gap between SFT(1,0.66) and SFT(1,0.51) narrows in Experiment 2. 
    In some cases (e.g., $n=45$ and $n=50$), SFT(1,0.51) achieves comparable objective values. The same happens for the SFT($n,\cdot$) algorithms. For instance, for $n \geq 40$, SFT ($n$,0.51) sometimes even provides better objective values than SFT ($n$,0.66).
    \item On Experiment Set 2, MS-OA continues to outperform RS-OA in terms of objective value. The objective values are slightly worse than in Experiment 1 but stay consistent overall.
    \item As for Experiment Set 1, the Insertion Heuristic still consistently provides the worst objective values, but maintaining by far the lowest computation time.
\end{itemize}
In general, computation times for all versions of the algorithms are significantly higher, especially for larger problem sizes. This is particularly noticeable for the SFT heuristics with $r=1$, while RS-OA and MS-OA show a smaller variation between the two experiments, with their computation times increasing more steadily. The objective values between the SFT algorithms are much closer in this second set of experiments, especially for the larger problem sizes. RS-OA and MS-OA still provide good quality solutions, with MS-OA generally outperforming RS-OA as in the first experiments. 

\begin{figure}[hbt]
     \centering
     \begin{tikzpicture}
 \begin{axis}
[
    title={\footnotesize Experiments Set 1: Objective-CPU times tradeoff when $n=65$},
    xlabel={\footnotesize Average CPU times},
    ylabel={\footnotesize Average objective value},
    legend pos=outer north east,
    ymajorgrids=true,
    grid style=dashed,
]
    \addplot [
        scatter,
        only marks,
        point meta=explicit symbolic,
        scatter/classes={
            Insert.H.={mark=square,blue},
            RS-OA={mark=triangle,red},
            MS-OA={mark=o,green},
            {SFT(1,0.66)}={mark=x,teal},
            {SFT(1,0.51)}={mark=otimes,teal},
            {SFT($n$,0.66)}={mark=diamond,teal},
            {SFT($n$,0.51)}={mark=oplus,teal},
            solver60={mark=Mercedes star flipped,brown},
            solver300={mark=Mercedes star,brown}
            },
     ] table [meta=label] {
    x      y       label
    6      3899    Insert.H.
    117    3596    RS-OA
    204    3580    MS-OA
    325    3625    {SFT(1,0.66)}
    118    3795    {SFT(1,0.51)}
    337    3624    {SFT($n$,0.66)}
    39     3875    {SFT($n$,0.51)}
    60     3767    solver60
    300    3700    solver300
    };
    \legend{Insert.H.,RS-OA,MS-OA,{SFT(1,0.66)},{SFT(1,0.51)},{SFT($n$,0.66)},{SFT($n$,0.51)},solver60,solver300}
 \end{axis}
 \end{tikzpicture}
    \caption{CPU time-Average objective values scatter plot for instances with $n=65$. 
    Symbols are Insertion Heuristic {\color{blue}$\square$},
            RS-OA {\color{red}$\triangle$},
            MS-OA {\color{green}$\circ$},
            SFT1-high {\color{teal}$\times$},
            SFT1-low  {\color{teal}$\otimes$},
            SFTmulti-high {\color{teal}$\diamond$},
            SFTmulti-low {\color{teal}$\oplus$},
            solver60 {\color{brown}$\Ydown$},
            solver300 {\color{brown}$\Yup$}.
    In this case, Insertion Heuristic, RS-OA and MS-OA constitute Pareto-efficient solution methods.
     }
     \label{fig:PO-plot}
 \end{figure}
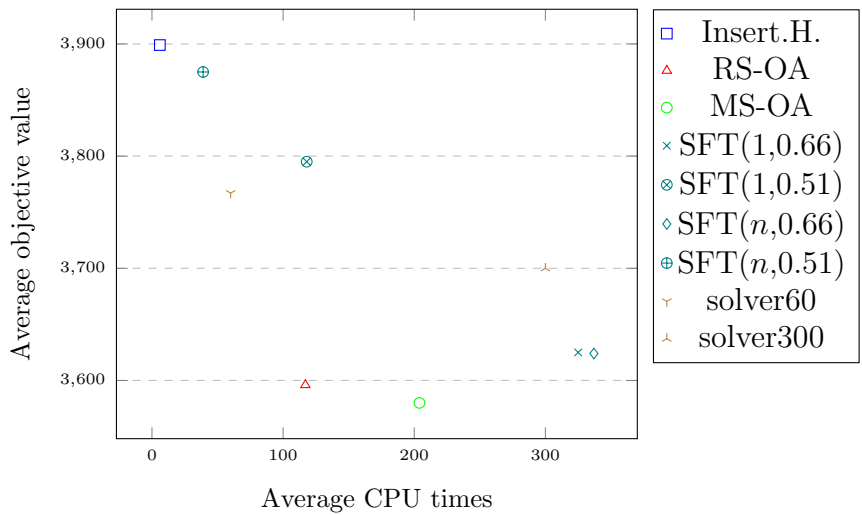

\medskip
A tradeoff analysis of the proposed methods can be summarized in the following points.  
For small problem sizes, SFT(1,0.66) or SFT($n$,0.66) provide the best objective values but at a higher computation cost. For moderate to large problem sizes MS-OA is the most balanced option, offering the best overall combination of speed and quality. RS-OA strikes a solid balance between speed and quality while for very fast approximate solutions the Insertion Heuristic is unmatched in speed but provides the least accurate solutions. This overall comparison suggests that MS-OA is the most versatile and well-rounded method across various problem sizes, particularly when both time and solution quality matter. Figure~\ref{fig:PO-plot} provides a visual illustration of the 
tradeoff between solution quality and computation time
for the instance class with \( n = 65 \) jobs of Experiment Set 1.

\section{Conclusions}\label{sec:conc}

In this work we investigated a 
 scheduling problem
encountered at a production line of an Austrian company manufacturing prefabricated house walls.
The problem can be modeled as a flexible permutation flow shop with blocking constraints, where at least one machine in the system has the ability to process operations from the preceding and succeeding stages. 
We identified some special cases that can be solved in polynomial 
time using graph-based algorithms.
To tackle the general problem, we proposed several ILP-based matheuristics and a constructive insertion heuristic. 
We evaluated the performance of the proposed heuristics through an extensive computational study based on realistic randomly generated instances 
providing insights into the effectiveness and efficiency of the proposed methods.
In particular, the computational tests indicate that, depending on the instance size and the time allowed for producing a solution, one approach may be preferable over another. 
Our constructive heuristic is by far the fastest of all tested algorithms and can even beat some of the MIP-based matheuristics w.r.t.\ solution quality.
Matheuristics based on fixing the job sequence first and then optimizing the assignment of shiftable operations yield good solutions within relatively low running times.
While there is no definitive winner, the proposed heuristic algorithms generally outperform the best-performing MIP models (and solvers) in terms of both solution quality and computation time.

As it is emphasized in Table~\ref{tab:complexity}, the complexity characterizations of two special cases of the addressed problem remain open, namely the cases of {\p} with $m=2$ and {\p} with $m=3$ and fixed job sequence.
In particular, when $m=2$, if one does not consider flexibility, an optimal solution can be found in polynomial time. On the other hand, if one has flexibile operation-to-machine assignment without the blocking constraint, the problem is NP-hard. 
Moreover, flexibility implies the NP-hardness of the flow shop problem with $m=2$ even if the sequence of jobs is fixed. However the same problem with blocking constraints is easy (as shown in Section~\ref{sec:fixedsorting}). We do not know yet if {\p} with $m\ge 3$ and fixed job sequence becomes hard.

Beside considering these open  questions concening computational complexity, future research might focus on the design of classical metaheuristic methods and study their performance in comparison to the approaches presented in this work.
 
\bibliographystyle{abbrv}
\bibliography{wall}

\end{document}